\documentclass[11pt]{article}

	\usepackage{amsmath, amssymb, amscd, latexsym, epsfig}
\newcommand{\arrow}{\rightarrow}

\newcommand{\cL}{{\cal L}}

\newcommand{\bb}{\mathbb}
\newcommand{\cx}{{\bb C}}

\newcommand{\integers}{{\bb Z}}

\newcommand{\natls}{{\bb N}}
\newcommand{\ratls}{{\bb Q}}
\newcommand{\reals}{{\bb R}}
\newcommand{\proj}{{\bb P}}
\newcommand{\makefig}[3]{
	\begin{figure}[htbp]
        \refstepcounter{figure}
	\label{#2}
        \begin{center}
		~#3~\\
		\medskip
                {\sf Figure \thefigure.  #1}
        \end{center}
	\end{figure}
}
\thicklines
\newcommand{\maketab}[3]{
	\begin{figure}[htbp]
        \refstepcounter{figure}
	\label{#2}
        \begin{center}
		#3~\\
		\bigskip
                {\sf Table \thefigure.  #1}
        \end{center}
	\end{figure}
}
\renewcommand{\bold}[1]{\smallskip \noindent {\bf \boldmath #1 }\nopagebreak[4]}

\newcommand{\qed}{\nopagebreak[4]\hfill
\rule{2mm}{2.5mm} \bigskip \pagebreak[2]}

\newcommand{\asyto}{\sim}

\newcommand{\bdry}{\partial}

\newcommand{\brackets}[1]{\langle #1 \rangle}

\newcommand{\closure}{\overline}
\newcommand{\compos}{\circ}

\newcommand{\del}{\partial}

\newcommand{\dirsum}{\oplus}

\newcommand{\disjunion}{\sqcup}

\newcommand{\includesin}{\hookrightarrow}

\newcommand{\isom}{\cong}

\newcommand{\lap}{\Delta}

\newcommand{\mem}{\in}

\newcommand{\st}{\: : \:}         

\newcommand{\Abar}{{\overline{A}}}

\newcommand{\Xbar}{{\overline{X}}}
\newcommand{\Sbar}{{\overline{S}}}
\newcommand{\Tbar}{{\overline{T}}}
\newcommand{\Ubar}{{\overline{U}}}
\newcommand{\Vbar}{{\overline{V}}}

\newcommand{\zbar}{{\overline{z}}}

\newcommand{\omegabar}{{\overline{\omega}}}

\newcommand{\chat}{{\widehat{\cx}}}

\newcommand{\area}{\operatorname{area}}

\newcommand{\Aut}{\operatorname{Aut}}

\newcommand{\dist}{\operatorname{dist}}

\newcommand{\logplus}{{\log^+}}

\renewcommand{\mod}{\operatorname{mod}}

\newcommand{\ord}{{\operatorname{ord}}}

\newcommand{\Poly}{\operatorname{Poly}}

\renewcommand{\Re}{\operatorname{Re}}

\newcommand{\val}{\operatorname{val}}

\newcommand{\zed}{\integers}

\newcommand{\cC}{{\cal C}}

\newcommand{\cF}{{\cal F}}

\newcommand{\cO}{{\cal O}}

\newcommand{\cT}{{\cal T}}

\newtheorem{theorem}{Theorem}[section]
\newtheorem{prop}[theorem]{Proposition}
\newtheorem{lemma}[theorem]{Lemma}
\newtheorem{cor}[theorem]{Corollary}

\makeatletter
\def\cleardoublepage{\clearpage\if@twoside \ifodd\c@page\else
    \thispagestyle{plain}\hbox{}\newpage\if@twocolumn\hbox{}\newpage\fi\fi\fi}

\def\ps@headings{\let\@mkboth\markboth
  \def\@oddfoot{}%
  \def\@evenfoot{}%
  \def\@evenhead{\small \sc\thepage\hfil\leftmark}
  \def\@oddhead{\small \sc \rightmark\hfil\thepage}
  \def\chaptermark##1{{
    \edef\@tempa{\ifnum \c@secnumdepth >\m@ne \@chapapp\ \thechapter. \fi}%
    \expandafter \markboth \expandafter{\@tempa ##1}{}}}%
  \def\schaptermark##1{\markboth {##1}{##1}}%
  \def\sectionmark##1{{
    \edef\@tempa{\ifnum \c@secnumdepth >\z@ \thesection. \fi}%
    \expandafter \markright \expandafter{\@tempa ##1}}}}

\def\thebibliography#1{\section*{References\@mkboth
 {References}{References}}\list
 {[\arabic{enumi}]}{\settowidth\labelwidth{[#1]}\leftmargin\labelwidth
 \advance\leftmargin\labelsep
 \usecounter{enumi}}
 \def\newblock{\hskip .11em plus .33em minus .07em}
 \sloppy\clubpenalty4000\widowpenalty4000
 \sfcode`\.=1000\relax}

\newif\if@restonecol
\def\theindex{\@restonecoltrue\if@twocolumn\@restonecolfalse\fi
\columnseprule \z@
\columnsep 35pt\twocolumn[\@makeschapterhead{Index}]
 \@mkboth{Index}{Index}\thispagestyle{plain}\parindent\z@
 \parskip\z@ plus .3pt\relax\let\item\@idxitem}
\def\@idxitem{\par\hangindent 40pt}

\def\endtheindex{\if@restonecol\onecolumn\else\clearpage\fi}

\def\footnoterule{\kern-3\p@ 
 \hrule width .4\columnwidth 
 \kern 2.6\p@} 
\@addtoreset{footnote}{chapter} 
\long\def\@makefntext#1{\parindent 1em\noindent 
 \hbox to 1.8em{\hss$^{\@thefnmark}$}#1}

\@addtoreset{equation}{section}

\renewcommand{\l@section}{\@dottedtocline{0}{1.5em}{2.3em}}
\renewcommand{\l@subsection}{\@dottedtocline{1}{3.8em}{3.2em}}
\renewcommand{\l@subsubsection}{\@dottedtocline{2}{7.0em}{4.1em}}

\makeatother

\input{xy} 
\xyoption{all}

\newcommand{\MPoly}{\operatorname{MPoly}}

\newcommand{\PT}{\proj \cT}
\renewcommand{\dist}{d}
\newcommand{\degf}{D}
\newcommand{\makefiglong}[3]{
        \begin{figure}[htbp]
        \refstepcounter{figure}
        \label{#2}
	\begin{center}
                ~#3~\\
	\end{center}
                \medskip
                {\sf Figure \thefigure.  #1}
        \end{figure}
}

\newcommand{\stroke}[2]{
   \multicolumn{12}{c}{\vspace{-2.5em}} \\
   \multicolumn{#1}{c}{~} &
   \multicolumn{#2}{c}{\hrulefill}  \\
   \multicolumn{12}{c}{\vspace{-1.7em}} \\
}

\title{ \vspace{-1in}
{\bf Trees and the dynamics of polynomials}
	\vspace{.2in}
}

\author{Laura G. DeMarco and Curtis T. McMullen\footnote{
Research supported in part by the NSF.}}

\begin{document}

\maketitle

\tableofcontents

\newpage

\begin{abstract}
In this paper we study branched coverings of metrized,
simplicial trees $F : T \arrow T$
which arise from polynomial maps $f : \cx \arrow \cx$
with disconnected Julia sets.
We show that the collection of all such trees, up to scale,
forms a contractible space $\PT_\degf$ compactifying the moduli space of 
polynomials of degree $\degf$;
that $F$ records the asymptotic behavior of the multipliers of $f$;
and that any meromorphic family of polynomials over $\Delta^*$
can be completed by a unique tree at its central fiber.
In the cubic case we give a combinatorial enumeration of the trees 
that arise, and show that $\PT_3$ is itself a tree.  
\end{abstract}

\renewcommand{\abstractname}{Resum\'e}
\begin{abstract}
Dans ce travail, nous \'etudions des rev\^etements ramifi\'es d'arbres
m\'etriques
simpliciaux $F: T \arrow T$ qui sont obtenus \`a partir d'applications
polyn\^omiales
$f : \cx \arrow \cx$ poss\'edant un ensemble de Julia non connexe.
Nous montrons que la collection de tous ces arbres, \`a un facteur
d'\'echelle pr\`es, forme un espace contractile $\PT_\degf$ qui compactifier
l'espace des modules des polyn\^omes de
degr\'e $\degf$.
Nous montrons aussi que $F$ 
enrigistre le comportement asymptotique 
des multiplicateurs de $f$ et  que tout
famille m\'eromorphe de polyn\^omes d\'efinis sur $\Delta^*$ peut
\^etre compl\'et\'ee par un unique arbre comme sa fibre
centrale. Dans le cas cubique, nous donnons une \'enum\'eration
combinatoire des arbres ainsi obtenus et montrons que $\PT_3$ est
lui-m\^eme un arbre.
\end{abstract}

\newpage



\thispagestyle{empty}


\section{Introduction}

The basin of infinity of a polynomial map $f : \cx \arrow \cx$ 
carries a natural foliation and a flat metric with singularities,
determined by the escape rate of orbits.
As $f$ diverges in the moduli space of polynomials, this Riemann surface collapses along its
foliation to yield a metrized simplicial tree $(T,\dist)$, with limiting dynamics $F : T \arrow T$.

In this paper we characterize the trees that arise as limits, and show they
provide a natural boundary $\PT_\degf$ compactifying 
the moduli space of polynomials of degree $\degf$. 
We show that $(T,\dist,F)$ records the limiting behavior of the multipliers of $f$ at its periodic
points, and that any degenerating analytic family of polynomials
$\{f_t(z) \st t \mem \Delta^* \}$ can be completed by a unique tree at
its central fiber.
Finally we show that in the cubic case, the boundary of moduli space $\PT_3$
is itself a tree; and for any $\degf$, $\PT_\degf$ is contractible.

The metrized trees $(T,\dist,F)$ provide a counterpart, in the setting of iterated rational
maps, to the $\reals$-trees that arise as limits of hyperbolic manifolds.

\bold{The quotient tree.}
Let $f : \cx \arrow \cx$ be a polynomial of degree $\degf \ge 2$.  The points $z \mem \cx$ with bounded orbits under $f$ form the compact {\em filled Julia set} 
$$K(f) = \{z: \sup_n|f^n(z)|<\infty\};$$
its complement, $\Omega(f) = \cx\setminus K(f)$, is the {\em basin of infinity}.
The {\em escape rate} $G : \cx \arrow [0,\infty)$ is defined by
\begin{displaymath}
	G(z) = \lim_{n\to \infty} \degf^{-n} \logplus |f^n(z)| ;
\end{displaymath}
it is the Green's function for $K(f)$ with a logarithmic pole at infinity.
The escape rate satisfies $G(f(z)) = \degf G(z)$, and thus it gives a semiconjugacy
from $f$ to the simple dynamical system $t \mapsto \degf t$ on $[0,\infty)$.

Now suppose that the Julia set $J(f) = \bdry K(f)$ is disconnected; equivalently, suppose
that at least one critical point of $f$ lies in the basin $\Omega(f)$.  
Then some fibers of $G$ are also disconnected, although for each $t>0$ the fiber 
$G^{-1}(t)$ has only finitely many components.

To record the combinatorial information of the dynamics of $f$ on $\Omega(f)$,
we form the {\em quotient tree} $\Tbar$ by
identifying points of $\cx$ that lie in the same connected
component of a level set of $G$.
The resulting space carries an induced dynamical system
\begin{displaymath}
	F : \Tbar \arrow \Tbar.
\end{displaymath}
The escape rate $G$ descends to give the {\em height function} $H$ on $\Tbar$,
yielding a commutative diagram
\begin{displaymath}
	\xymatrix{\cx \ar[d]_{\pi} \ar[dr]^{G} & \\
			\Tbar \ar[r]_-H  &  [0,\infty)	}
\end{displaymath}
respecting the dynamics.  See Figure \ref{fig:crit} for an illustration of the trees in two examples.  Note that only a finite subtree of $\Tbar$ has been drawn in each case.  

The {\em Julia set} of $F:\Tbar\to\Tbar$ is defined by
\begin{displaymath}
	J(F) = \pi(J(f))= H^{-1}(0) . 
\end{displaymath}
It is a Cantor set with one point for each connected component of $J(f)$.
With respect to the measure of maximal entropy $\mu_f$, 
the quotient map $\pi : J(f) \arrow J(F)$ is almost injective, in the 
sense that $\mu_f$-almost every component of $J(f)$ is a single
point (Theorem \ref{thm:pointcomponents}).  In particular, there is no loss of information
when passing to the quotient dynamical system:

\begin{theorem}  
Let $f$ be a polynomial of degree $\degf\geq 2$ with disconnected Julia set.
The measure-theoretic entropy of $(f,J(f),\mu_f)$ and its quotient
$(F,J(F),\pi_*(\mu_f))$ are the same --- they are both $\log \degf$.
\end{theorem}

\makefiglong{Critical level sets of $G(z)$ for $f(z)=z^2+2$ and $f(z) = z^3 + 2.3z^2$, and 
the corresponding combinatorial trees.  The edges mapping with degree $>1$ are indicated.  The Julia set of $z^2+2$ is a Cantor set, while the Julia set of $z^3 + 2.3 z^2$ contains countably many Jordan curves.}
{fig:crit}
{ ~\hspace{-0.4in}\psfig{figure=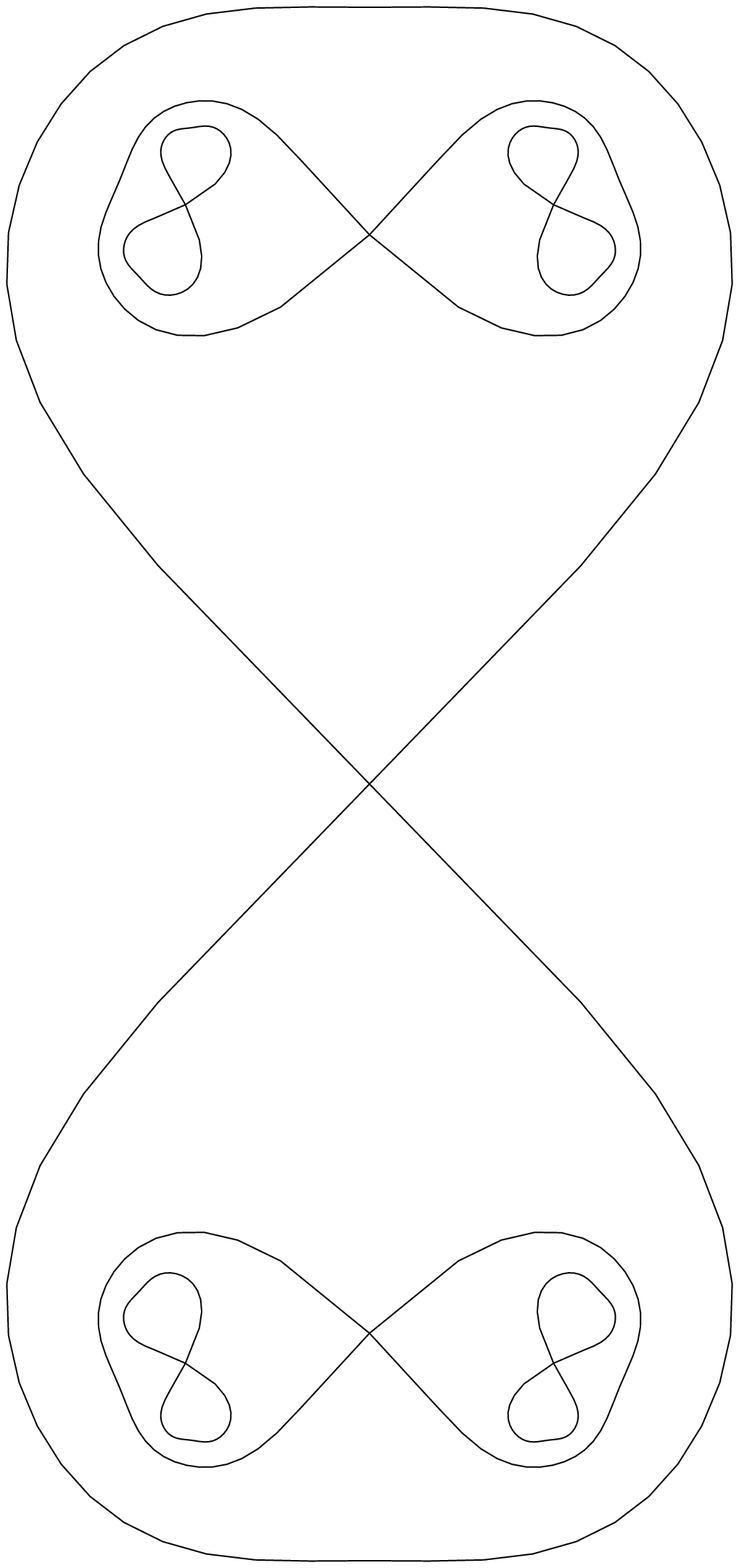,height=1.8in}
\hspace{0.5in}  \psfig{figure=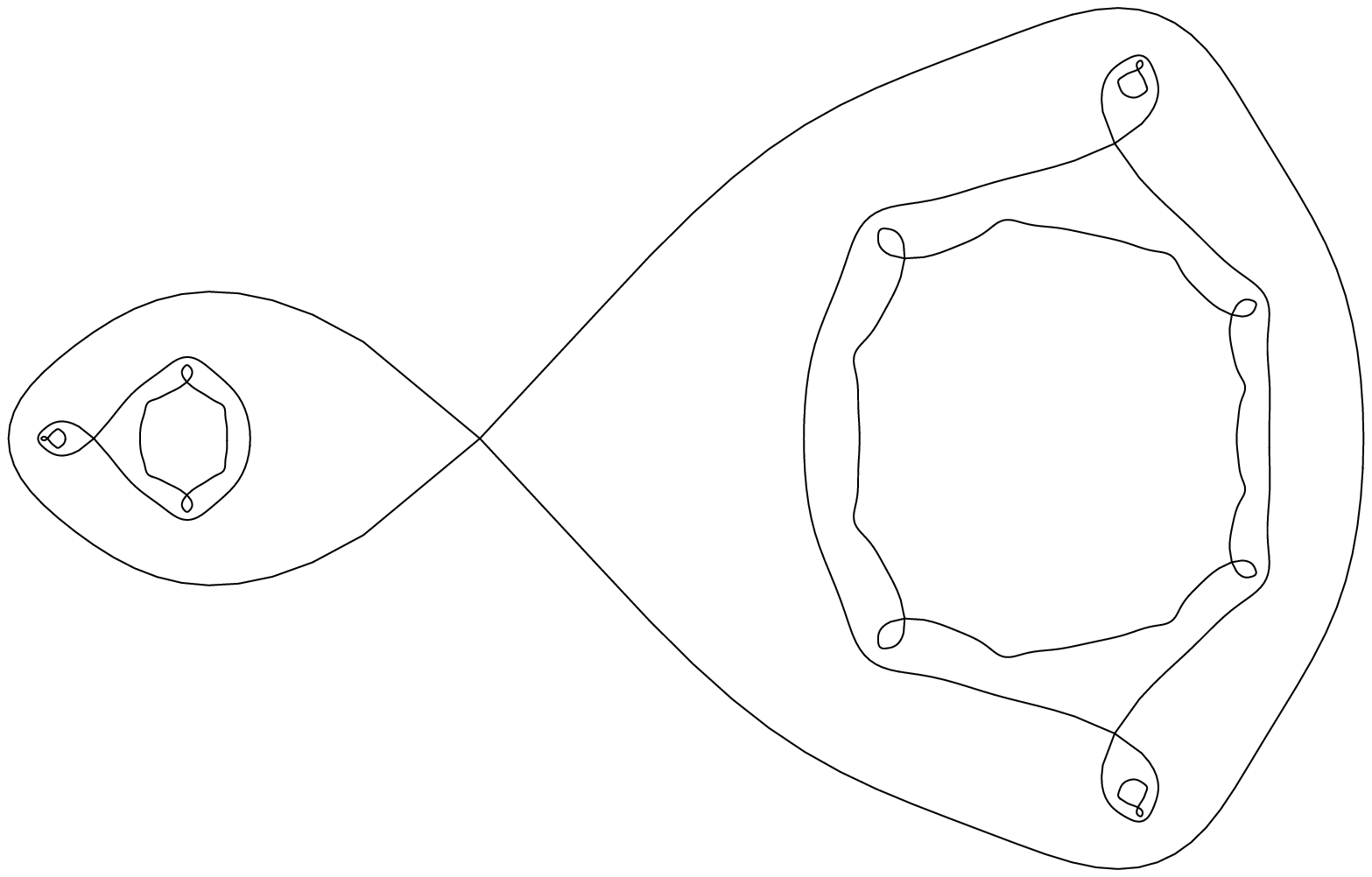,height=1.8in}\\
\bigskip \psfig{figure=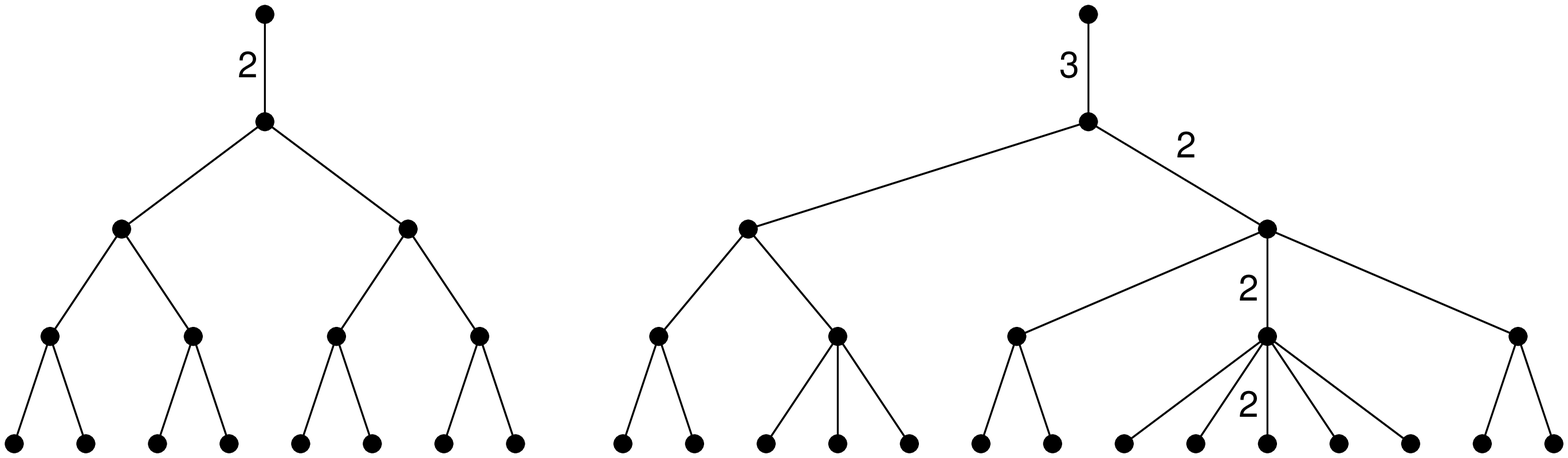,width=\textwidth}  }

The quotient of the basin of infinity,
\begin{displaymath}
	T = \pi(\Omega(f)) = H^{-1}(0,\infty),
\end{displaymath}
is an open subset of $\Tbar$ homeomorphic to a simplicial tree.  In fact $T$ carries
a canonical simplicial structure, determined by the conditions:
\begin{enumerate}
	\item
$F : T \arrow T$ is a simplicial map,
	\item
the vertices of $T$ consist of the grand orbits of the branch points of $T$, and
	\item
the height function $H$ is linear on each edge of $T$.
\end{enumerate}
Details can be found in \S\ref{sec:trees}.

The {\em height metric} $\dist$ is a path metric on $T$,
defined so that adjacent vertices of $T$
satisfy
$$	\dist(v,v') = |H(v)-H(v')| . $$
We refer to the triple
\begin{displaymath}
	\tau(f) = (T,\dist,F)
\end{displaymath}
as the {\em metrized polynomial-like tree} obtained as the quotient of $f$.
The space $\Tbar$ is the metric completion of $(T,\dist)$. 

In \S\ref{sec:trees}, we introduce an {\em abstract} metrized polynomial-like tree with dynamics 
$F : T \arrow T$, and in \S\ref{sec:realization} we show:

\begin{theorem}  
\label{thm:Tsurjective}
Every metrized polynomial-like tree $(T,\dist, F)$
arises as the quotient $\tau(f)$ of a polynomial $f$.
\end{theorem}

\noindent
Special cases of Theorem \ref{thm:Tsurjective} were proved by Emerson; Theorems 9.4 and 10.1 of \cite{Emerson:trees} show that any tree with just one escaping critical point (though possibly of high multiplicity) and with divergent sums of moduli can be realized by a polynomial.

\bold{Spaces of trees and polynomials.}
Let $\cT_\degf$ denote the space of isometry classes of metrized 
polynomial-like trees $(T,\dist,F)$ of degree $\degf$.
The space $\cT_\degf$ carries a natural geometric topology, defined by convergence of
finite subtrees.
There is a continuous action of $\reals_+$ on $\cT_\degf$ by rescaling the metric $\dist$,
yielding as quotient the projective space
\begin{displaymath}
	\PT_\degf = \cT_\degf/\reals_+ .
\end{displaymath}
In \S\ref{sec:topology} we show:

\begin{theorem} 
\label{thm:compactcontractible}
The space $\PT_\degf$ is compact and contractible.
\end{theorem}

Now let $\MPoly_\degf$ denote the moduli space of polynomials of degree $\degf$, the space
of polynomials modulo conjugation by the affine automorphisms of $\cx$.  The conjugacy class
of a polynomial $f$ will be denoted $[f]$.  
The space $\MPoly_\degf$ is a complex orbifold, finitely covered by $\cx^{\degf-1}$.
The {\em maximal escape rate}
\begin{displaymath}
	M(f) = \max\{G(c): f'(c)=0\}
\end{displaymath}
depends only the conjugacy class of $f$; Branner and Hubbard observed in \cite{Branner:Hubbard:cubicsI}
that $M$ descends to a continuous and proper map
$M : \MPoly_\degf \arrow [0,\infty)$.  

The {\em connectedness locus} $\cC_\degf \subset \MPoly_\degf$ is the subset
of polynomials with connected Julia set; it coincides with the locus $M(f)=0$ and is therefore
a compact subset of $\MPoly_\degf$.
We denote its complement by
\begin{displaymath}
	\MPoly_\degf^* = \MPoly_\degf \setminus \cC_\degf.
\end{displaymath}
The metrized polynomial-like tree $\tau(f) = (T,\dist, F)$ depends only on 
the conjugacy class of $f$, so $\tau$ induces a map 
  $$\tau:  \MPoly_\degf^* \arrow \cT_\degf.$$ 
Note that the compactness of the connectedness locus $\cC_\degf$ implies
that every divergent sequence in $\MPoly_\degf$ will eventually
be contained in the domain of $\tau$. 

There is a natural action of $\reals_+$ on $\MPoly_\degf^*$ obtained by
`stretching' the complex structure on the basin of infinity. 
In \S\ref{sec:continuity} and \S\ref{sec:compact} we show:

\begin{theorem}
\label{thm:quotmap}
The map $\tau: \MPoly_\degf^* \arrow \cT_\degf$ is continuous, proper, surjective, 
and equivariant with respect to the action of $\reals_+$ by stretching of
polynomials and by metric rescaling of trees.
\end{theorem}

\begin{theorem}
\label{thm:boundary}
The moduli space of polynomials admits a natural compactification
$\closure{\MPoly}_\degf = \MPoly_\degf \cup \,\PT_\degf $
such that 
\begin{itemize}
\item  $\MPoly_\degf$ is dense in $\closure{\MPoly}_\degf$, and 
\item  the iteration map $[f] \mapsto [f^n]$ extends 
	continuously to $\closure{\MPoly}_\degf \arrow \closure{\MPoly}_{\degf^n}$.
\end{itemize}
\end{theorem}

\bold{Periodic points.}  Fix a polynomial $f$ with disconnected 
Julia set, and let $\tau(f) = (T,\dist,F)$ be its metrized polynomial-like 
tree.
The {\em modulus metric} $\delta$ is another useful path metric 
on $T$, defined on adjacent vertices by
	$$\delta(v,v') = 2\pi\mod(A)$$
where $A = \pi^{-1}(e) \subset \Omega(f)$
is the annulus lying over the open edge $e$ joining $v$ to $v'$.
(Here $\mod(A)=h/c$ when $A$ is conformally 
a right cylinder of height $h$ and circumference $c$.)
Let $p \mem J(F)$ be a fixed point of $F^n$.
The {\em translation length} of $F^n$ at $p$ is defined by
\begin{displaymath}
	L(p,F^n) = \lim_{v \arrow p} \delta(v,F^n(v)) ,
\end{displaymath}
where the limit is taken over vertices $v \mem T$ along the 
unique path from $\infty$ to $p$.
In \S\ref{sec:multipliers} we establish:

\begin{theorem}  
\label{thm:multiplierbound}
Let $z \mem \cx$ be a fixed point of $f^n$, and let $p=\pi(z) \mem J(F)$.
Then the log-multiplier of $z$ and translation length at $p$ satisfy
\begin{displaymath}
        L(p,F^n) \le \log^+ |(f^n)'(z)| \le L(p,F^n) + C(n,D),
\end{displaymath}
where $C(n,D)$ is a constant depending only on $n$ and $D$.
\end{theorem}
The argument also shows the periodic points $p \mem J(F)$ with 
$L(p,F^n)>0$ are in bijective correspondence with the 
periodic points $z \mem J(f)$ that form singleton components of the
Julia set (see Proposition \ref{prop:repelling}).
In particular, we have the following curious consequence: 

\begin{cor}
All singleton periodic points in $J(f)$ are repelling.
\end{cor}

A metrized tree $(T,\dist,F)$ is {\em normalized} if the distance from the highest
branched point $v_0\in T$ to the Julia set $J(F)$ is 1.    
In \S\ref{sec:compact}, we introduce a notion of {\em pointed} convergence of 
polynomials and trees, and we use Theorem \ref{thm:multiplierbound} to prove:

\begin{theorem}
\label{thm:onemultiplier}
Suppose $[f_k]$ is a sequence in $\MPoly_\degf$ 
which converges to the normalized tree $(T,\dist,F)$ in the boundary $\PT_\degf$.
Let $z_k\in\cx$ be a sequence of fixed-points of $(f_k)^n$ converging to $p \mem \Tbar$.
Then the translation length of $F^n$ at $p$ is given by
\begin{displaymath}
	L(p,F^n) = \lim_{k\to\infty} \frac{\log^+|(f_k^n)'(z_k)|}{M(f_k)} \cdot
\end{displaymath}
\end{theorem}
Recall that $M$ is the maximal escape rate, and so $M(f_k)$ tends 
to infinity as $k\to\infty$.  

\bold{Metrizing the basin of infinity.}
The holomorphic 1-form $\omega = 2 \del G$ provides a dynamically
determined conformal metric $|\omega|$ on the basin of infinity $\Omega(f)$, with singularities at
the escaping critical points and their inverse images.  
In this metric $f$ is locally expanding by a factor of $\degf$,
and a neighborhood of infinity is isometric to a cylinder 
$S^1 \times [0,\infty)$ of radius one.

Let $c(f)$ denote (one of) the fastest escaping critical point(s) of $f$,
so that $M(f) = G(c(f))$.
Let $\Xbar(f)$ denote the metric completion of $(\Omega(f),|\omega|)$,
rescaled so the distance from $c(f)$ to the boundary $\Xbar(f)\setminus 
\Omega(f)$ is $1$.
In \S\ref{sec:GH} we show:

\begin{theorem}  \label{thm:GHboundary}
If $[f_k]$ converges to the normalized tree $(T,\dist,F)$ in $\closure{\MPoly}_\degf$,
then
\begin{displaymath}
	(\Xbar(f_n),c(f_n))  \arrow (\Tbar,v_0)
\end{displaymath}
in the Gromov-Hausdorff topology on pointed metric spaces, where $v_0$ is the 
highest branch point of the tree $\Tbar$, 
and $f_n : \Xbar(f_n) \arrow \Xbar(f_n)$ converges to $F : \Tbar \arrow \Tbar$.
\end{theorem}

\bold{Algebraic limits.}
Theorems \ref{thm:onemultiplier} and  \ref{thm:GHboundary}
show the space of trees 
$\PT_\degf = \bdry \MPoly_\degf$
is large enough to
record the growth of multipliers at periodic points and the limiting geometry of
the basin of infinity.  The next result shows it is small enough that
any holomorphic map of the punctured disk
\begin{displaymath}
	\Delta^* \arrow \MPoly_\degf
\end{displaymath}
which is meromorphic at $t=0$ extends to a continuous map $\Delta \arrow \closure{\MPoly}_\degf$
(see \S\ref{sec:families}).

\begin{theorem}
\label{thm:algebraic}
Let $f_t(z) = z^\degf + a_2(t) z^{\degf-2} + \cdots + a_\degf(t)$
be a holomorphic family of polynomials over $\Delta^*$, whose coefficients
have poles of finite order at $t=0$.
Then either:
\begin{itemize}
	\item
$f_t(z)$ extends holomorphically to $t=0$, or 
	\item
the conjugacy classes of $f_t$ in $\MPoly_\degf$ converge
to a unique normalized tree $(T,\dist,F) \mem \PT_\degf$ as $t \arrow 0$.
\end{itemize}
In the latter case the edges of $T$ have rational length, and hence
the translation lengths of its periodic points are also rational.
\end{theorem}

The rationality of translation lengths is related to valuations:  it reflects the fact 
that the multiplier $\lambda(t)$ of a periodic point of $f_t$ is given by a Puiseux series 
\begin{displaymath}
	\lambda(t) = t^{p/q} + [\text{higher order terms}]
\end{displaymath}
at $t=0$; compare \cite{Kiwi:trees}.

\bold{Cubic polynomials.}
We conclude by examining the topology of $\PT_\degf$ for $\degf=3$.
Given a partition $1 \le \sum p_i \le \degf-1$, let
\begin{displaymath}
	\cT_\degf(p_1,\ldots,p_N) \subset \cT_\degf
\end{displaymath}
denote the locus where the escaping critical points fall into $N$ grand orbits,
each containing $p_i$ points (counted with multiplicity).
Each connected component of
$\proj \cT_\degf(p_1,\ldots,p_N)$ is an open simplex of dimension $N-1$ 
(see Proposition \ref{prop:htmetrics}),
and each component of
$$\bigcup_{\sum p_i=e} \PT_\degf(p_1,\ldots,p_N)$$ 
is a simplicial complex of dimension $e-1$.
In \S\ref{sec:cubics} we show:

\begin{theorem}
The boundary $\PT_3$ of the moduli space of cubic polynomials
is the union of an infinite simplicial tree $\PT_3(2) \cup \PT_3(1,1)$ and its 
set of ends $\PT_3(1)$.
\end{theorem}

\noindent
We also give an algorithm for constructing $\proj\cT_3$ 
via a combinatorial encoding of its vertices $\PT_3(2)$.  

\bold{Notes and references.}
Branner and Hubbard initiated the study
of the tree-like combinatorics of Julia sets, and its ramifications for the moduli space of polynomials, especially cubics,
using the language of tableaux
\cite{Branner:Hubbard:cubicsI}, \cite{Branner:Hubbard:cubicsII}; see also \cite{Milnor:local:connectivity}.

Trees for polynomials of the type we consider here 
were studied independently by 
Emerson \cite{Emerson:trees}; see also \cite{Emerson:harmonic}.
Other connections between trees and complex dynamics appear in
\cite{Shishikura:rings} and \cite{Przytycki:Skrzypczak}.

Trees also arise naturally from limits of group actions on hyperbolic spaces
(see e.g.  \cite{Morgan:Shalen:I}, \cite{Morgan:Shalen:survey}, \cite{Morgan:icm86}, \cite{Otal:book:fibered}, \cite{Paulin:rtrees}).
For hyperbolic surface groups, the space of limiting $\reals$-trees coincides with Thurston's 
boundary to Teichm\"uller space, and the
translation lengths on the $\reals$-tree record the limiting behavior of 
the lengths of closed geodesics.
These results motivated the formulation of 
Theorems \ref{thm:onemultiplier} and \ref{thm:multboundary}.   
The theory of $\reals$-trees can be developed for
limits of proper holomorphic maps the unit disk as well \cite{McMullen:rtrees}.
For a survey of connections between rational maps and Kleinian groups, see \cite{McMullen:curdev}.

We would like to thank the referee for useful comments.  


\section{Abstract trees with dynamics}
\label{sec:trees}

In this section we discuss polynomial-like maps $F$ on simplicial trees.
We show $F$ has a naturally defined set of critical points,
a canonical invariant measure $\mu_F$ on its Julia set $J(F)$,
and a finite-dimensional space of compatible metrics.

In \S\ref{sec:realization} we will show that
every polynomial-like tree actually comes from a polynomial.

\bold{Trees.}
A {\em simplicial tree} $T$ is a nonempty, connected, locally finite, 1-dimensional simplicial complex 
without cycles.
The set of vertices of $T$ will be denoted by $V(T)$, and the set of
(unoriented, closed) edges by $E(T)$.
The edges adjacent to a given vertex $v \mem V(T)$ form a finite set
$E_v(T)$, whose cardinality $\val(v)$ is the {\em valence} of $v$.

The space of {\em ends} of $T$, denoted $\bdry T$, is the compact, 
totally disconnected space
obtained as the inverse limit of the set of connected components 
of $T-K$ as $K$ ranges over all finite subtrees.
The union $T \cup \bdry T$,
with its natural topology, is compact.

\bold{Branched covers.}
We say a map $F : T_1 \arrow T_2$ between simplicial trees is a
{\em branched covering} if:
\begin{enumerate}
	\item
$F$ is proper, open and continuous; and
	\item
$F$ is simplicial (every edge maps linearly to another edge).
\end{enumerate}
These conditions imply that $F$ extends continuously to the boundary of $T$, yielding an
open, surjective map
\begin{displaymath}
	F : \bdry T_1 \arrow \bdry T_2 .
\end{displaymath}

\bold{Local and global degree.}
A {\em local degree function} for a branched covering $F$ is a map
\begin{displaymath}
	\deg : E(T_1) \cup V(T_1) \arrow \{1,2,3\ldots\}
\end{displaymath}
satisfying, for every $v \mem V(T_1)$, the inequality
\begin{equation}   \label{eq:degineq}
	2 \deg(v) - 2 \geq \sum_{e \mem E_v(T_1)} (\deg(e)-1) ,
\end{equation}
as well as the equality
\begin{equation}
\label{eq:degeq}
	\deg(v) = \sum_{e' \mem E_v(T_1) : F(e')=F(e)} \deg(e') 
\end{equation}
for every $e \mem E_v(T_1)$.
These conditions imply that $(F,\deg)$ is locally modeled on a $\deg(v)$
branched covering map between spheres.  The tree maps arising from
polynomials always have this property (\S\ref{sec:polys}).

In terms of the local degree, the global degree $\deg(F)$ is defined by:
\begin{equation}
\label{eq:degF}
	\deg(F) = \sum_{F(e_1)=e_2} \deg(e_1) = \sum_{F(v_1)=v_2} \deg(v_1)
\end{equation}
for any edge $e_2$ or vertex $v_2$ in $T_2$.
It is easy to verify that this expression is independent
of the choice of $e_2$ or $v_2$, using (\ref{eq:degeq}) and
connectedness of $T_2$.

\bold{Polynomial-like tree maps.}
Now let $F : T \arrow T$ be the dynamical system given by a branched covering map of
a simplicial tree to itself.
Two points of $T$ are in the same {\em grand orbit} if $F^n(x)=F^m(y)$ for some $n,m>0$.

We say $F$ is {\em polynomial-like} if:
\begin{enumerate}
	\item[I.]
There is a unique isolated end $\infty \mem \bdry T$;
	\item[II.]
There exists a local degree function compatible with $F$; 
	\item[III.]
The tree $T$ has no endpoints (vertices of valence one); and
	\item[IV.]
The grand orbit of any vertex includes a vertex of valence $3$ or more.
\end{enumerate}
We will later see that the local degree function is unique 
(Theorem \ref{thm:deg}).

Here are some basic properties of a polynomial-like $F:T\to T$ that
follow quickly from the definitions.
\begin{enumerate}
	\item
We have $\val(v) \ge \val(F(v)) \ge 2$ for all $v \mem V(T)$.

	\item
If $\val(v)=2$, then its adjacent edges satisfy
\begin{equation}
\label{eq:e1e2}
	\deg(e_1)=\deg(e_2)=\deg(v).
\end{equation}

	\item
If $\val(v)\ge 3$, then $F$ is a local homeomorphism at $v$
if and only if $\deg(v)=1$.

	\item
The {\em Julia set}
\begin{displaymath}
	J(F) = \bdry T - \{\infty\}
\end{displaymath}
is homeomorphic to a Cantor set.  Note that $J(F)$ is 
compact, totally disconnected and perfect, and $J(F)$ is nonempty because
$T$ has no endpoints.

	\item
Every point of $J(F)$ is a limit of vertices of valence three or more.

	\item
The extended map $F:\Tbar\to\Tbar$ is finite-to-one.  This follows from (\ref{eq:degF}).

	\item
The point $\infty \mem \bdry T$ is totally invariant;
that is, $F^{-1}\{\infty\} = \{\infty\}$.  This follows from 
the fact that $\infty$ is the unique isolated point in $\bdry T$ and 
$F|\bdry T$ is finite, continuous and surjective.
\end{enumerate}

\bold{Combinatorial height.}
Since the end $\infty \mem \bdry T$ is isolated, every vertex close enough to
$\infty$ has valence two.
The {\em base} $v_0 \mem T$ is the unique vertex of
valence 3 or more that is closest to $\infty \mem \bdry T$.
This vertex splits $T$ into a pair of subtrees
\begin{displaymath}
	T = (J(F),v_0] \cup [v_0,\infty)
\end{displaymath}
meeting only at $v_0$.   The subtree
$[v_0,\infty)$ is an infinite path converging to $\infty \mem \bdry T$;
the remainder $(J(F),v_0]$  accumulates on the Julia set.

The {\em combinatorial height function}
\begin{displaymath}
	h : V(T) \arrow \zed
\end{displaymath}
is defined by $|h(v)| = $ the minimal number of edges
needed to connect $v$ to $v_0$; its sign is determined by the condition that
$h(v) \ge 0$ on $[v_0,\infty)$ while $h(v) \le 0$ on $(J(F),v_0]$.
(Equivalently, $-h(v)$ is a horofunction measuring the 
number of edges between $v$ and $\infty$, normalized so $h(v_0)=0$.)

\begin{lemma}
\label{lemma:h}
There is an integer $N(F)>0$ such that the combinatorial height satisfies
\begin{equation}
\label{eq:h}
	h(F(v)) = h(v) + N(F).
\end{equation}
\end{lemma}

\bold{Proof.} 
Since every internal vertex of $(v_0,\infty)$ has degree two,
$F|(v_0,\infty)$ is a homeomorphism.  Since $F(\infty)=\infty$,
we must have $F(v_0,\infty) \subset (v_0,\infty)$.
(The image cannot contain $v_0$ since $\val(v_0) \ge 3$.)
Consequently (\ref{eq:h}) holds for all $v$ with $h(v) \ge 1$, 
with $N(F) \ge 0$.  In fact we have $N(F)>0$;
otherwise any vertex close enough to $\infty$ would be
totally invariant (since $\infty$ is), contradicting our assumption
that its grand orbit contains a vertex of valence 3 or more.

To see (\ref{eq:h}) holds globally, just note that
$h(F(v))$ can have no local maximum, by openness of $F$.
\qed

\begin{cor}
We have $F^n(x) \arrow \infty$ for every $x \mem T$,
and the quotient space $T/\brackets{F}$
is a simplicial circle with $N(F)$ vertices.
\end{cor}
(To form the quotient, we identify every grand orbit to a single point.)

\bold{Local degree.}
Every vertex $v \mem V(T)$ has a unique {\em upper edge}, leading
towards $\infty$; and one or more {\em lower edges}, leading to
vertices of lower height.

\begin{lemma}
\label{thm:e}
The upper edge $e$ of any vertex $v$ satisfies
$\deg(e)=\deg(v)$; and $\deg(v)=\deg(F)$ when $h(v) \ge 0$.
\end{lemma}

\bold{Proof.} 
By (\ref{eq:h}), only the upper edge $e$ at $v$ can
map to the upper edge at $F(v)$, and thus
$\deg(e)=\deg(v)$.
For the second statement, suppose $h(v) \ge 0$ and
$F(v)=F(v')$; then $h(v')=h(v)$, so $v'=v$.
Applying equation (\ref{eq:degF}) with $e_2=F(e)$,
we obtain $\deg(F) = \deg(v)$.
\qed

\begin{cor}
The degree function is increasing along any sequence of consecutive
vertices converging to $\infty$.
\end{cor}

\bold{Critical points.}
We define the {\em critical multiplicity} of a vertex $v \mem V(T)$ by
\begin{displaymath}
	m(v) = 2\deg(v)-2 - \sum_{e\in E_v(T)} (\deg(e)-1),
\end{displaymath}
which is non-negative by (\ref{eq:degineq}).  
Similarly, if $v_i$ is a sequence of consecutive vertices converging to a
point $p \mem J(F)$, then $\deg(v_i)$ is decreasing and we define
\begin{displaymath}
	m(p) = \lim \deg(v_i) - 1 \ge 0.
\end{displaymath}
If $x \mem V(T) \cup J(F)$ and $m(x)>0$, we say 
$x$ is a {\em critical point} of multiplicity $m(x)$.

\begin{lemma}
The base $v_0\in T$ is a critical point.
\end{lemma}

\bold{Proof.}   
The base $v_0$ has degree $\deg(F)$ and 
lower edges $e_1, \ldots, e_k$ with $k>1$ and 
$F(e_1) = \cdots = F(e_k)$.  The 
critical multiplicity is therefore $m(v_0) = k-1>0$.
\qed

\begin{lemma}  \label{criticalpoints}
The total number of critical points of $F$, counted with multiplicity,
is $\deg(F)-1$.
The degree of an edge is one more than the number of critical points below it,
counted with multiplicities.
\end{lemma}

\bold{Proof.} 
Using Lemma \ref{thm:e}, the critical multiplicity can be computed as
$$m(v) = \deg(v) - 1 - \sum_{E_l(v)} (\deg(e)-1),$$
where $E_l(v)$ is the collection of lower edges of $v$.  Furthermore, if 
$e_u$ is the upper edge of $v$, then 
$$\deg(e_u) = \deg(v) = m(v) + 1 + \sum_{E_l(v)} (\deg(e)-1).$$
For each lower edge, we can replace $\deg(e)$ with a similar expression
involving the critical multiplicity of the vertex below it and degrees of its
lower edges.  Continuing inductively, we conclude that the degree of $e_u$
is exactly one more than the number of critical points below it.  In particular, 
since $\deg(v) = \deg(F)$ for all vertices above the base $v_0$, 
the total number of critical points is $\deg(F)-1$.
\qed

\begin{cor}
The edges with $\deg(e)>1$ form the convex hull of the critical points
union $\{\infty\}$.
\end{cor}

\begin{cor}
We have $N(F) \le \deg(F)-1$.
\end{cor}

\bold{Proof.} 
Every vertex of valence two is in the forward orbit
of a vertex of valence three or more, and hence in the forward orbit of a critical point.
\qed

\bold{Uniqueness of the local degree.}
Let $J(F,v)$ denote the subset of the Julia set lying below a given vertex $v \mem V(T)$.  That is, $J(F,v)$ is the collection of all ends $p\in J(F)$ such that the unique path joining $\infty$ and $p$ passes through $v$.  Because $F$ takes the lower edges of a vertex surjectively to the lower edges of the image vertex, we have 
  $$F(J(F,v)) = J(F,F(v)).$$
We can now show:

\begin{theorem}
\label{thm:deg}
If $F : T \arrow T$ is polynomial-like, then its 
local degree function is unique.  The degree $\deg(v)$ of 
a vertex $v$ is the topological degree of $F|J(F,v)$, counting critical points
with multiplicity.
\end{theorem}

\bold{Proof.} 
Fix a vertex $v$ and an end $q$ below $F(v)$.  Set $w_0 = F(v)$, and let $w_i$ denote the consecutive sequence of vertices tending to $q$ with combinatorial height $h(w_i) = h(w_0) -i$.  Let $e$ be the upper edge of $w_1$ (so it is a lower edge of $F(v)$).   By (\ref{eq:degeq}), we have
  $$\deg(v) = \sum_{e'\in E_v(T):  F(e') = e} \deg(e').$$
From Lemma \ref{thm:e}, $\deg(e') = \deg(v')$ where $e'$ is the upper edge to vertex $v'$, and consequently,
  $$\deg(v) = \sum_{v' \mbox{ below } v, \; F(v') = w_1} \deg(v').$$
Proceeding inductively on the combinatorial height, we have 
  $$\deg(v) = \sum_{v' \mbox{ below } v, \; F(v') = w_i} \deg(v')$$
for every $i\geq 1$.  Passing to the limit, we see that $\deg(v)$ records the number of preimages of the end $q$, counted with multiplicities.  Finally, this implies uniqueness of the degree function because there are only finitely many critical points.  The degree of $v$ is the number of preimages in $J(F,v)$ of a generic point in $J(F,F(v))$.
\qed

\bold{Vertex counts.}
Because of the preceding result, 
(\ref{eq:degF}) gives an unambiguous definition of the
global degree of a polynomial-like $F$.
Next we show $T$ has controlled exponentially growth below its base.
Let $V_k(T) = \{v \mem V(T) \st h(v) = -k N(F)\}$.   

\begin{lemma} \label{thm:counts}  
Let $\degf = \deg(F)$.  For any $k \ge 0$, we have:
\begin{displaymath}
	\degf^k \ge |V_k(T)| \ge 2 + \degf + \degf^2 + \cdots + \degf^{k-1} \ge \frac{\degf^k}{\degf-1} \cdot
\end{displaymath}
\end{lemma}

\bold{Proof.} 
The upper bound follows from the fact that $|V_0(T)|=1$ and
$|V_{k+1}(T)| \le \degf |V_k(T)|$, since $F(V_{k+1}(T)) = V_{k}(T)$.
Since there are at most $(\degf-2)$ critical points below the base of the tree,
and $\deg(v)=1$ unless there is a critical point at or below $v$, we also have:  
\begin{displaymath}
	|V_{k+1}(T)| \ge \degf |V_k(T)| - (\degf-2),
\end{displaymath}
which gives the lower bound.
\qed

\bold{Invariant measure.}
The {\em mass function} $\mu : V(T) \arrow \ratls$
is characterized by the conditions
\begin{equation}
\label{eq:r}
	\mu(F(v)) = \frac{\deg(F)}{\deg(v)} \cdot \mu(v) 
\end{equation}
for all $v \mem V(T)$, and $\mu(v) = 1$ when $h(v) \ge 0$.
These properties determine $\mu(v)$ uniquely,
since the forward orbit of every vertex converges to $\infty$.

Recall that $J(F,v)$ denotes the subset of the Julia set
lying below a given vertex $v \mem V(T)$.
Note that $F|J(F,v)$ is injective if $\deg(v)=1$.

By induction on the combinatorial height,
one can readily verify that if 
$v_1,\ldots,v_s$ are the vertices immediately below $v$, then
\begin{displaymath}
	\mu(v) = \mu(v_1) + \mu(v_2) + \cdots \mu(v_s) .
\end{displaymath}
Consequently, there is a unique Borel probability measure
$\mu_F$ on $J(T)$ satisfying
\begin{displaymath}
	\mu_F(J(F,v)) = \mu(v)
\end{displaymath}
for all $v \mem V(T)$.

\begin{lemma}
The probability measure $\mu_F$ is invariant under $F$.
\end{lemma}

\bold{Proof.} 
By (\ref{eq:degF}), if $F^{-1}(v) = \{v_1,\ldots,v_s\}$, then
$\deg(F) = \sum_1^s \deg(v_i)$ and thus
\begin{displaymath}
	\mu(v) = \sum_1^s \deg(v_i)/\deg(F) = \sum_1^s \mu(v_i) .
\end{displaymath}
Consequently we have
\begin{displaymath}
	\mu_F(F^{-1}(J(F,v)) =
	\sum_1^s \mu_F(J(F,v_i)) = \sum_1^s \mu(v_i) = \mu(v) =
		\mu_F(J(F,v)).
\end{displaymath}
Since open sets of the form $J(F,v)$ generate the Borel algebra of $J(F)$,
$\mu_F$ is invariant.
\qed

The exponential growth of $T$ gives an {\em a priori} diffusion
to the mass of $\mu_F$.

\begin{lemma}
\label{lem:diffuse}
For any vertex $v \mem V_k(T)$, we have:
\begin{displaymath}
	\mu_F(J(F,v)) \le \left( \frac{\degf-1}{\degf} \right)^k .
\end{displaymath}
\end{lemma}

\bold{Proof.} 
The vertex $v$ maps in $k$ iterates to $v_0$,
which satisfies $\mu(v_0)=1$.  Along the way the,
the degree is bounded by $(\degf-1)$, and thus $\mu(v) 
\le ((\degf-1)/\degf)^k$ by (\ref{eq:r}).
\qed

\begin{cor}
The measure $\mu_F$ has no atoms.
\end{cor}

\begin{cor}
For any Borel set where $F|A$ is injective, we have:
\begin{equation}
\label{eq:A}
	\mu_F(F(A)) = \deg(F) \cdot \mu_F(A).
\end{equation}
\end{cor}

\bold{Proof.} 
By (\ref{eq:r}) this Corollary holds when $A=J(F,v)$ and
$\deg(v)=1$; and since there are only finitely many critical points
in $J(F)$,
$\deg(v) = 1$ for some vertex above almost any point in $J(F)$.   
\qed

\bold{Univalent maps.}
The degree function for the map $F^n : T \arrow T$ is given by
\begin{displaymath}
	\deg(v,F^n) = \deg(v) \cdot \deg(F(v)) \cdots \deg(F^{n-1}(v)) .
\end{displaymath}
We say $F^n$ is {\em univalent} at $v$ if
$\deg(v,F^n) = 1$.  The next result shows `almost every' vertex can be
mapped univalently up to a definite height.

\begin{lemma}
\label{thm:univ}
For almost every $x \mem J(F)$ there exists a $k \ge 0$ such that
for each $i \ge 0$, the map
$F^i$ is univalent at the vertex $v \mem V_{k+i}(T)$ lying above $x$.
\end{lemma}

\bold{Proof.} 
Let $C_i$ denote the union of $J(F,v)$ for all vertices $v \mem V_i(T)$
lying above critical points of $F$.
Since there are most $(\degf-2)$ critical points below the base of $T$, 
Lemma \ref{lem:diffuse} gives
\begin{displaymath}
	\mu_F(C_i) \leq (\degf -2) ((\degf-1)/\degf)^i .
\end{displaymath}
It is easily verified that the lemma holds for all $x$ in
\begin{displaymath}
	J_k = J(F) - \bigcup_{i=1}^\infty F^{-i}(C_{k+i}).
\end{displaymath}
Since $F$ is measure preserving and $\sum \mu_F(C_i) < \infty$, we have
$\mu_F(J_k) \arrow 1$ as $k \arrow \infty$, and thus 
the lemma holds for almost every $x \mem J(F)$.
\qed

\bold{Entropy and ergodicity.}
We can now show:

\begin{theorem}
\label{thm:entropy}
The invariant measure $\mu_F$ for $F|J(F)$ is ergodic,
and its entropy is $\log \deg(F)$.
\end{theorem}

\bold{Proof.} 
Let $A \subset J(F)$ be an $F$-invariant Borel set of positive measure.
Then the density of $A$ in $J(F,v)$ tends to 1 as $v$ approaches
almost any point of $A$.  Pick a point of density $x \mem A$ where
Lemma \ref{thm:univ} also holds.  
Then there exists a vertex $v \mem V_k(T)$ such that arbitrarily small neighborhoods
of $x$ map univalently onto $J(F,v)$.
By (\ref{eq:A}) the density of $A$ is preserved under univalent
maps, and hence there $A$ has density $1$ in some $J(F,v)$.
But $F^k(J(F,v)) = J(F)$, and thus $A$ has full measure.

Similarly, Lemma \ref{thm:univ} implies
that for almost every $x \mem J(F)$,
the vertices $v_n \mem V_n(T)$ lying above $x$ satisfy, for $n \ge 0$,
\begin{displaymath}
	- \log \mu(v_n) = n \log \deg(F) + O(1) .
\end{displaymath}
The entropy of $F$ is thus $\log \deg(F)$ by the Shannon-McMillan-Breiman theorem
\cite{Parry:book:entropy}.
\qed

See \cite{Brolin:measure}, \cite{FLM:entropy}, 
\cite{Lyubich:entropy} and \cite{Gromov:entropy} for analogous results for 
polynomials and rational maps.

\bold{Height metric.}
A {\em path metric} $d(x,y)$ on a simplicial tree $T$ is a
metric satisfying
\begin{displaymath}
	d(x,y)+d(y,z)=d(x,z) 
\end{displaymath}
whenever $y$ lies on the unique arc connecting $x$ to $z$ 
(cf. \cite[\S 1.7]{Gromov:book:var}).
We will also require that a path metric is linear on edges (with respect to
the simplicial structure).  Then $d$ is determined by the lengths
$d(e) = d(x,y)$ it assigns to edges $e = [x,y] \mem E(T)$.

A {\em height metric} $\dist$ for $(T,F)$ is a path metric satisfying
\begin{displaymath}
	\dist(F(e)) = \deg(F) \cdot \dist(e).
\end{displaymath}
A height metric is uniquely determined by the lengths it assigns
to the edges
\begin{displaymath}
	e_i  = [v_{i-1},v_i],\;\;i=1,2,\ldots,N(F)	
\end{displaymath}
joining the consecutive vertices $v_0, \ldots, v_{N(F)} = F^{N(F)}(v_0)$,
since this list includes exactly one edge from
each grand orbit in $E(T)$.  The lengths of these edges can
be arbitrary, and therefore:

\begin{prop}
\label{prop:htmetrics}
The set of height metrics $\dist$ compatible with $(T,F)$ is parameterized by
$\reals_+^{N(F)}$.  
\end{prop}

Since any path leading to the Julia set has length bounded by
$O(\sum \degf^{-n})$, the space
\begin{displaymath}
	\Tbar = T \cup J(F)
\end{displaymath}
is homeomorphic to the metric completion of $(T,\dist)$.
Moreover, the {\em height function} $H : \Tbar \arrow [0,\infty)$,
defined by
\begin{displaymath}
	H(x) = \dist(x,J(F)),
\end{displaymath}
satisfies $H(F(x)) = \deg(F) \cdot H(x)$.

\bold{Modulus metric.}
A height metric determines a unique {\em modulus metric} $\delta$,
characterized the conditions
\begin{displaymath}
	\delta(F(e)) = \deg(e) \cdot \delta(e)
\end{displaymath}
for all $e \mem E(T)$, and
by $\delta(e) = \dist(e)$ for edges $e$ in $[v_0,\infty)$.
Note that if $e$ is the upper edge at $v$, we have $\delta(e) \mu(v) = \dist(e)$.

\begin{lemma}
\label{thm:inf}
Almost every $x \mem J(F)$ lies at infinite distance from $v_0$ in
the modulus metric.
\end{lemma}

\bold{Proof.} 
By Lemma \ref{thm:univ}, for almost every $x$ there is a $k \ge 0$ and
a sequence of consecutive vertices $v_i \arrow x$, each of which can be mapped univalently
up to a vertex in $V_k(T)$.  Since $V_k(T)$ is finite, the correspond upper edges $e_i$ have
$\delta(e_i)$ bounded below, and thus $\sum \delta(e_i) = \infty$.
\qed

\bold{Summary.}
For later applications we will focus on the metric space $(T,\dist)$ and
its dynamics $F$.  Since the vertices of $T$
are the grand orbits of its branch points, the simplicial structure
and its further consequences are already implicit in this data.

\begin{theorem}
The metric space $(T,\dist)$ and the continuous map $F : T \arrow T$
uniquely determine:
\begin{enumerate}
	\item
the simplicial structure of $T$, 
	\item
the degree function on its vertices and edges,
	\item
the set of critical points with multiplicities,
	\item
the height function $H : T \arrow \reals$,
	\item
the modulus metric $\delta$ and
	\item
the invariant measure $\mu_F$ on $J(F)$.
\end{enumerate}
\end{theorem}
We refer to the triple $(T,\dist,F)$ as a {\em metrized polynomial-like tree}.


\section{Trees from polynomials}
\label{sec:polys}

In this section we discuss the relationship between a polynomial $f(z)$ and the
quotient dynamical system $\tau(f) = (T,\dist,F)$.

\bold{Foliations, metrics and measures.}
Let $f : \cx \arrow \cx$ be a polynomial of degree $\degf \ge 2$, with escape
rate $G(z) = \lim \degf^{-n} \logplus |f^n(z)|$ as in the Introduction.

The level sets of $G$ determine a foliation $\cF$ of the basin of infinity
$\Omega(f)$, with transverse invariant measure $|dG|$.
The holomorphic 1-form 
\begin{displaymath}
	\omega = 2\del G \asyto dz/z
\end{displaymath}
determines a flat metric $|\omega|$ making the leaves of $\cF$ into closed geodesics.
The distribution
\begin{displaymath}
	\mu_f =  (2\pi)^{-1} \lap G
\end{displaymath}
gives the harmonic measure on the Julia set $J(f)$, as well as the
probability measure of maximal entropy, $\log \degf$ 
\cite{Lyubich:entropy}.

The length of a closed leaf $L$ of $\cF$ determines the
measure of the Julia set inside the disk $U$ it bounds; namely, we have:
\begin{equation}
\label{eq:L}
	2\pi \mu_f(U) = \int_U \lap G = \int_L |\omega|  
\end{equation}
by Stokes' theorem.
The foliation and metric have isolated singularities along the
grand orbits of the critical points in $\Omega(f)$.

\bold{The quotient tree.}
As in the Introduction,  let $T$ be the space obtained by collapsing each 
leaf of $\cF$ to a single point, and let
\begin{displaymath}
	\pi : \Omega(f) \arrow T 
\end{displaymath}
be the quotient map.  We make $T$ into a metric space by defining
\begin{displaymath}
	\dist(\pi(a),\pi(b)) = \inf\int_a^b |\omega|,
\end{displaymath}
where the infimum is over all paths joining $a$ to $b$.
Since $f$ preserves the level sets of $G$, it descends to give
a map $F : T \arrow T$.

\begin{theorem}
If $\Omega(f)$ contains a critical point, then
$(T,\dist,F)$ is a metrized polynomial-like tree.
\end{theorem}

\bold{Proof.} 
Since the map $G : \Omega(f) \arrow (0,\infty)$ is proper, with a discrete
set of critical points, the quotient $T$ is a tree.
Its branch points come from the critical points of $G$,
which coincide with the backwards orbits of critical point of $f$ in $\Omega(f)$.
The maximum principle implies $T$ has no endpoints.

Since $f|\Omega(f)$ is open and proper, so is $F|T$.
The projections of the grand orbits of the critical points determine
a discrete set of vertices $V(T)$, giving $T$ a compatible simplicial
structure.  The level set of $G$ near $z=\infty$ are connected,
so $z=\infty$ gives an isolated end of $T$.
On the other hand, the Julia set $J(f)$ is contained in the closure of the
grand orbit of any critical point in $\Omega(f)$, so the remaining ends
of $T$ are not isolated.

Finally we show $F$ has a compatible degree function.
Since $G$ is a submersion over $T-V(T)$,
the preimage of the midpoint of an edge $e$ is a smooth loop
$L(e) \subset \Omega(f)$.  Given a vertex $v$,
let $S(v) \subset \Omega(f)$ denote the 
compact region bounded by the loops $L(e)$ for edges adjacent to $v$.
Note that $f : S(v) \arrow S(F(v))$ is a branched covering map,
with branch points only in the interior. 
Defining
\begin{displaymath}
	\deg(e) = \deg(f | L(e)),
	\;\;\;\;
	\deg(v) = \deg(f | S(v)),
\end{displaymath}
we see the degree axioms (\ref{eq:degeq}) and (\ref{eq:degineq})
follow from the Riemann-Hurwitz formula and the fact that
$\deg(f|S(v)) = \deg(f|\bdry S(v))$.
\qed

\bold{Dictionary.}
Recall from \S\ref{sec:trees} that $(T,\dist,F)$ determines a
set of critical points, a height function, a modulus metric and an invariant measure.
These objects correspond to $f$ as follows.

\begin{enumerate}
	\item
The critical vertices of $T$ are the images of the critical points of $f$.
Every vertex lies in the grand orbit of a critical point.
	\item
The height function $H : T \arrow (0,\infty)$ satisfies $H(\pi(x)) = G(x)$
as in the Introduction.
	\item
The preimage of the interior of $e \mem E(T)$ is an open annulus
$A(e)$ foliated by smooth level sets of $G$.  In the $|\omega|$-metric,
this annulus has height $\dist(e)$ and satisfies
\begin{equation}
\label{eq:modA}
	2\pi \mod(A) = \delta(e).
\end{equation}
	\item
The degree of an edge $e$ is the same as the degree of $f : A(e) \arrow A(F(e))$.
	\item
If $e$ is the upper edge of $v$, then the circumference of $A(e)$ is
given by $(2\pi) \mu_F(v)$.
	\item
The quotient map $\pi$ extends continuously to a map $\pi : \cx \arrow \Tbar$
sending $K(f)$ to $J(F)$ by collapsing its components to distinct, single points.
By the preceding observation and (\ref{eq:L}), this map satisfies
\begin{equation}
\label{eq:push}
	\pi_*(\mu_f) = \mu_F .
\end{equation}
	\item
The measures $\mu_f$ and $\mu_F$ have the same entropy, namely
$\log \degf$.
	\item
The critical points in $J(F)$ are the images of the critical points in $K(f)$.
\end{enumerate}

\bold{Functoriality.}
We remark that the tree construction is functorial:  a conformal 
conjugacy from $f(z)$ to $g(z)$ determines an isometry
between the quotient trees $\tau(f)$ and $\tau(g)$, respecting the dynamics.
Similarly, if $\tau(f) = (T,\dist,F)$ then $\tau(f^n) = (T_n,\dist_n,F_n)$ is naturally isometric 
to $(T,\dist,F^n)$.

\bold{Singletons.}
We say $x \mem J(f)$ is a {\em singleton} if
$\{x\}$ is a connected component of $J(f)$.

\begin{theorem} \label{thm:pointcomponents}
If $J(f)$ is disconnected, then $\mu_f$-almost every point $x \mem J(f)$ 
is a singleton.
\end{theorem}

\bold{Proof.} 
Let $x \mem J(f)$ and $y = \pi(x) \mem J(F)$.
By Theorem \ref{thm:inf} and (\ref{eq:push}), $y$ is almost
surely at infinite distance from $v_0$ in the modulus metric.
This means there is a sequence of consecutive edges $e_i$ leading
to $y$ with $\sum \delta(e_i) = \infty$.
Thus by (\ref{eq:modA}), the disjoint annuli $A(e_i) \subset \Omega(f)$ nested
around $x$ satisfy $\sum \mod(A_i) = \infty$, and 
therefore $x$ is a singleton.
\qed

\begin{cor}
The map $\pi : (J(f),\mu_f) \arrow (J(F),\mu_F)$ becomes a bijection after
excluding sets of measure zero.
\end{cor}

This gives another proof that $\pi$ preserves measure-theoretic entropy.

\bold{Remark.}
Qiu and Yin and, independently, Kozlovski and van Strien have recently shown that for any polynomial $f(z)$,
all but countably many components of $J(f)$ are singletons \cite{Qiu:Yin}, \cite{Kozlovski:vanStrien}.
For a rational map, however, the Julia set can be
homeomorphic to the product of a Cantor set with a circle, 
as for $f(z) = z^2 + \epsilon/z^3$ with $\epsilon$ small
\cite{McMullen:aut}.
Another proof of Theorem \ref{thm:pointcomponents}, using \cite{Qiu:Yin},
appears in \cite{Emerson:harmonic}.


\section{Multipliers and translation lengths}
\label{sec:multipliers}

Let $(T,\dist,F)$ be the quotient tree of a polynomial $f(z)$.
In this section we introduce the translation lengths $L(p,F^n)$,
and establish:

\begin{theorem}
\label{thm:mult}
Let $z \mem \cx$ be a fixed point of $f^n$, and let $p=\pi(z) \mem J(F)$.
Then the log-multiplier of $z$ and translation length at $p$ satisfy
\begin{displaymath}
        L(p,F^n) \le \log^+ |(f^n)'(z)| \le L(p,F^n) + C(n,D),
\end{displaymath}
where $C(n,D)$ is a constant depending only on $n$ and $D$.
\end{theorem}
This result is a restatement of Theorem \ref{thm:multiplierbound}.

We remark that the inequality $\log^+ |(f^n)'(z)| \ge L(p,F^n)$
follows easily from subadditivity of the modulus,
using the fact that a path in $T$ corresponds to a sequence of 
nested annuli in $\cx$.
For the reverse inequality, we must show these annuli
are glued together efficiently. 

\bold{Definitions.}
Let $f(z)$ be a polynomial with disconnected Julia set.  The 
{\em log-multiplier} of a periodic point $z$ of period $n$ is the 
quantity $\log^+|(f^n)'(z)|$.  

Let $(T,\dist,F)$ be the quotient tree of $f$.
Let $p \mem J(F)$ be a fixed point of $F^n$
and $v_i$ be a sequence of consecutive vertices converging to $p$.
Using the modulus metric (see (\ref{eq:modA})), we define the {\em translation length} of $F^n$ at $p$ by:
\begin{displaymath}
	L(p,F^n) = \lim_{i\to\infty} \delta(v_i,F^n(v_i)) .
\end{displaymath}
If the forward orbit of $p$ contains a critical point of $F$, then
$L(p,F^n)=0$.
Otherwise, $F^n$ is univalent at $v_i$ for all $i$ sufficiently large,
and hence it eventually acts by an isometric translation on
the infinite path leading to $p$.
In this case we say $p$ is a {\em repelling periodic point.}
We have $L(p,F^n) > 0$ since every point in $T$ converges to
infinity under iteration.

\begin{prop}
\label{prop:repelling}
The repelling periodic points in $J(F)$ correspond bijectively
to the singleton repelling periodic points in $J(f)$.
\end{prop}

\bold{Proof.} 
If $p=\pi(z)$ is a repelling periodic point then the path
from $v_0$ to $\pi(z)$ has infinite length in the modulus metric,
so $z$ is a singleton.  Any edge $e$ sufficiently close to 
$p$, along the path from $p$ to $\infty$, 
gives a nested pair of annuli encircling $p$ and mapping by degree one:
\begin{displaymath}
	A(e) \stackrel{f^n}{\to} A(F^n(e));
\end{displaymath}
thus $|(f^n)'(z)| > 1$ by the Schwarz lemma.

Conversely, if $z \mem J(f)$ is a periodic
singleton then it cannot be a critical point of $f^n$, so
$p=\pi(z)$ is repelling.
\qed

\bold{Polynomial-like maps.}
A proper holomorphic map $f : U_1 \arrow U_0$
between regions in the plane is {\em polynomial-like} if
$\closure{U_1}$ is a compact subset of $U_0$ and
$U_0-\Ubar_1$ is an annulus.
To begin the proof of Theorem \ref{thm:mult}, we show:

\begin{theorem}
\label{thm:polylike}
Let $f : U_1 \arrow U_0$ be a polynomial-like map of degree $d \ge 2$,
whose critical values lie in $U_1$.  Let $U_2 = f^{-1}(U_1)$, and suppose
\begin{displaymath}
	1/m < \mod(U_0-\Ubar_1) < \mod(U_0-\Ubar_2) < m.
\end{displaymath}
Then the fixed points of
$f$ satisfy
\begin{displaymath}
	|f'(p)| \le C(d,m) .
\end{displaymath}
\end{theorem}

\bold{Proof.} 
By the Riemann mapping theorem we can assume $U_0$ is the unit
disk $\Delta$ and $p=0$.
We can then write
\begin{displaymath}
	f = B \compos h
\end{displaymath}
where $h : U_1 \arrow U_0$ is
degree one, $B : U_0 \arrow U_0$ is degree $d$,
and $B(0)=h(0)=p$.
Let $V = B^{-1}(U_1)$, so $U_2 = h^{-1}(V)$.
Then we have 
\begin{displaymath}
	\mod(U_1-\Ubar_2) = \mod(U_0-\Vbar) = (1/d) \mod(U_0-\Ubar_1) \ge 1/(dm),
\end{displaymath}
since $h : (U_1 - \Ubar_2) \arrow (U_0 - \Vbar)$ is an isomorphism,
and $B : (U_0-\Vbar) \arrow (U_0-\Ubar_1)$ is a covering map of degree $d$.

Since $0 \mem U_2$ and $\mod(U_0-\Ubar_2)  = \mod(\Delta-\Ubar_2) \le m$,  
there is a point $q \mem U_2$ with $|q| > r(m,d) > 0$. 
Since the annulus $U_1-\Ubar_2$ has modulus $\ge 1/(dm)$ and encloses
$\{p,q\} = \{0,q\}$, the region $U_1$ contains a ball of radius 
$r'(d,m) = C(m)r(m,d) > 0$ about $p=0$.
Finally, since $h$ maps $U_1$ into $\Delta$, the Schwarz lemma implies
\begin{displaymath}
	|f'(0)| = |h'(0)| \cdot |B'(0)|
		\le |h'(0)| \le 1/r'(m,d),
\end{displaymath}
as required.
\qed

\bold{Counterexample.}
We emphasize that the preceding result is false if we only require
$1/m < \mod(U_0-\Ubar_1) < m$.

To see this, let $B : \Delta \arrow \Delta$ be a fixed degree two Blaschke
product with $B(0)=0$ and with its unique critical value at $z=-1/6$.
Let $M_r(z) = (z+r)/(1+rz)$, and let $A_r(z) = (z-r)/3$, where $0<r<1$.
Then $A_r(M_r(\Delta)) = U_r$ is the disk of radius $1/3$ centered at
$-r/3$, so it contains the critical value of $B$.  Moreover,
\begin{displaymath}
	h_r = (A_r \compos M_r)^{-1} : U_r \arrow \Delta
\end{displaymath}
is a degree one map, with $h_r(0)=0$ and $h_r'(0) = 3/(1-r^2)$.
Thus 
\begin{displaymath}
	f_r = B \compos h_r : U_r \arrow \Delta
\end{displaymath}
is a polynomial-like map of degree 2,
with critical values in $U_r$ and with $\mod(\Delta-\Ubar_r)$ bounded
above and below.  On the other hand $f_r(0)=0$, and the multiplier
\begin{displaymath}
	|f'_r(0)| = 3|B'(0)|/(1-r^2) 
\end{displaymath}
tends to infinity as $r \arrow 1$.

\bold{Consecutive annuli.}
Next we give an estimate for the modulus of an
annulus $A \subset \cx$ formed from 
consecutive annuli $A_1,\ldots,A_n$ of the kind
that arise from the tree construction.

\makefig{A nest of consecutive annuli.}{fig:nest}{
\psfig{figure=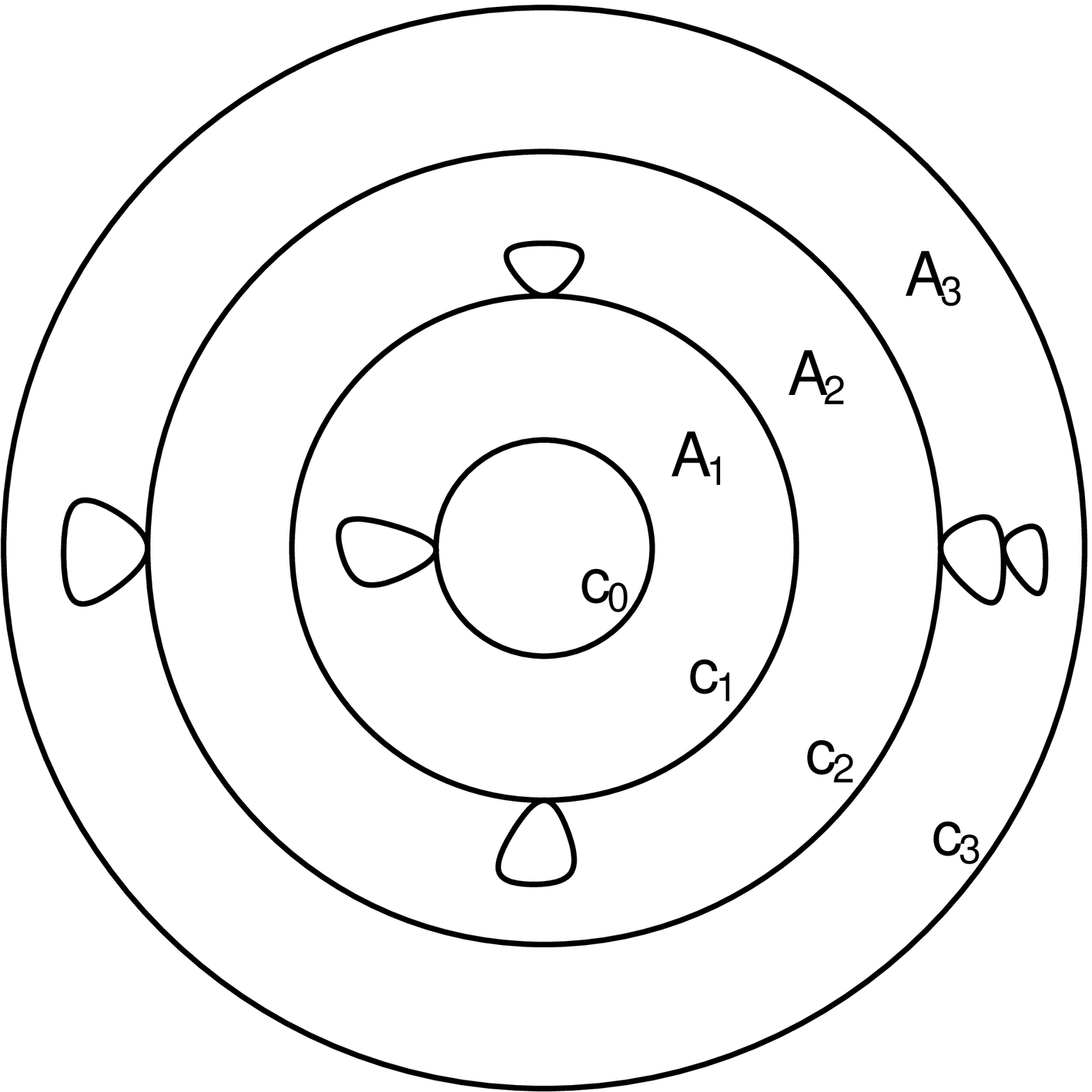,height=2.2in}}

Let $\bigcup_1^n A_i \subset A \subset \cx$ be a set
of disjoint nested annuli $A_i$ inside an annulus $A$.
Assume:
\begin{enumerate}
	\item
Each annulus has piecewise smooth inner and
outer boundaries, $\bdry_-A_i$ and $\bdry_+ A_i$,
	\item
The outer boundary of $A_i$ is a Jordan curve,
made up of finitely many segments of the inner boundary of $A_{i+1}$
(so long as $i<n$); 
	\item
There is a continuous conformal metric
$\rho = \rho(z) |dz|$ on $\bigcup_1^n \Abar_i$,
making each annulus $A_i$ into a flat right cylinder
of height $h_i$ and circumference $c_i$; 
	\item
The boundary of $A$ is a pair Jordan curves, with
$\bdry_-A \subset \bdry_-A_1$ and
$\bdry_+A = \bdry_+ A_n$.
\end{enumerate}
These conditions imply $\mod(A_i)=h_i/c_i$.
Letting $c_0$ denote the $\rho$-length of $\bdry_-A$,
we have 
\begin{displaymath}
	c_0 \le c_1 \le \cdots \le c_n.
\end{displaymath}

\begin{theorem}
\label{thm:nest}
The modulus of $A$ satisfies:
\begin{equation}
\label{eq:summod}
	\sum_1^n \mod(A_i) \le \mod(A) \le
	3n(c_n/c_0)^2 +
	\sum_1^n \mod(A_i)  .
\end{equation}
\end{theorem}

\bold{Proof.} 
The first inequality is standard; for the second,
we will use the method of extremal length 
(cf. \cite{Lehto:Virtanen:book}).

Let us say an annulus $A_i$ is {\em short} if 
$h_i < 2c_0+c_i$; otherwise it is {\em tall}.  
Define a conformal metric $\sigma$ on $\Abar$
by setting $\sigma = (1/c_0) \rho$ on all the
short annuli, and on cylindrical collars 
of $\rho$-height $c_0$ at the two ends of the tall annuli.  
Between the collars of each tall annulus $A_i$,
let $\sigma = (1/c_i)\rho$.
Extend $\sigma$ to the rest of $A$ by setting it equal to zero.

Let $\Gamma$ denote the set of all rectifiable
loops in $A$ separating its boundary components.
It is now straightforward to verify that 
\begin{displaymath}
	L_\sigma(\gamma) = \int_\gamma \sigma \ge 1
\end{displaymath}
for all $\gamma \mem \Gamma$.

To see this, first suppose $\gamma$ meets the
region between the collars of a tall annulus $A_i$.
If $\gamma$ is contained in $A_i$ then it must separate the
boundary components of $A_i$, so $L_\rho(\gamma) \ge c_i$,
and thus $L_\sigma(\gamma) \ge 1$ 
(since $1/c_0 > 1/c_i$).
Otherwise $\gamma$ must cross one of the collars of $A_i$;
but each collar has $\sigma$-height one, so 
again $L_\sigma(\gamma) \ge 1$.

Now suppose $\gamma \cap \bigcup \Abar_i$
is covered by short annuli and the collars of 
tall annuli.
On this region $\sigma = (1/c_0) \rho$.
Consider the foliation $\cF$ of $\bigcup \Abar_i$
by geodesics in the flat $\rho$-metric, which
start at $\bdry_- A$ and proceed perpendicular to the boundary
in the outward direction.
Any $\gamma \mem \Gamma$ must cross all the leaves of $\cF$.  
By construction the leaves are parallel, with constant separation,
within the short annuli and collars of tall annuli.
Thus the projection of $\gamma \cap \cF$
(along leaves of $\cF$) to $\bdry_-A$ is 
$\sigma$-distance decreasing, and thus
\begin{displaymath}
	L_\sigma(\gamma) \ge L_\sigma(\bdry_-A) = (1/c_0)c_0 = 1
\end{displaymath}
in this case as well.

Since the modulus of $A$ is the reciprocal of the
extremal length of $\Gamma$, we have:
\begin{displaymath}
	\mod(A) = 1/\lambda(\Gamma) \le
	\left( \int_A \sigma^2 \right)\left/
	\left( \inf_\Gamma \int_\gamma \sigma \right)^2 \right.
	\le \sum \area_\sigma(A_i).
\end{displaymath}
Each short annulus has height $h_i \le 2c_0+c_i \le 3c_n$,
so it contributes area $h_ic_i/c_0^2 \le 3(c_n/c_0)^2$.
Each tall annulus contributes area at most
$h_i/c_i + 2c_0c_i/c_0^2 \le \mod(A_i) + 3(c_n/c_0)^2$;
summing over $i$, we obtain (\ref{eq:summod}).
\qed

\bold{Torus shape.}
Let $f : \bdry_-A \arrow \bdry_+A$ be a piecewise smooth 
homeomorphism preserving orientation, and expanding the metric 
$\rho$ linearly by a factor of $c_n/c_0$.
Let 
\begin{displaymath}
	T = \Abar/f
\end{displaymath}
be the complex torus obtained by
gluing together corresponding points, and let
$B \subset T$ be an annulus of maximum modulus homotopic to $A$.

A straightforward modification of the proof above yields:

\begin{theorem}
\label{thm:torus}
We have $\sum \mod(A_i) \le \mod(B)$.
In addition, we have
\begin{displaymath}
	\mod(B) \le 
	3n(c_n/c_0)^2  + \sum \mod(A_i) 
\end{displaymath}
provided $\mod(A_1) \ge 3$.
\end{theorem}
The condition on $\mod(A_1)$ implies that $A_1$ is a tall annulus, and
hence it cannot be crossed by loops with $L_\sigma(\gamma) \le 1$.

\bold{Bounds on multipliers.}
We can now complete the proof of Theorem \ref{thm:mult}.
It suffices to treat the case where $z$ is a fixed point of $f$.

\begin{lemma}  \label{lemma:onesided}
We have $L(p,F) \le \log^+ |f'(z)|$.
\end{lemma}

\bold{Proof.} 
The statement is clear if $L(p,F)=0$.  Otherwise both $p$ and $z$ are 
repelling fixed points (by Proposition \ref{prop:repelling}),
and the degree of $F$ is one near $p$.
Let $e_i$, $i \mem \zed$, be the unique path of consecutive
edges in $T$ connecting $p$ to $\infty$; it satisfies
\begin{displaymath}
	F(e_i) = e_{i+n},
\end{displaymath}
where $n = N(F) \le \degf-1$.
Note that $\deg(e_i)$ is monotone increasing, and equal to $1$
for all $i$ sufficiently small.
After shifting indices we can assume $\deg(e_n)=1$;
then $\deg(e_i)=1$ for all $i \le n$, and we have
\begin{displaymath}
	L(p,F) = \sum_1^n \delta(e_i) .
\end{displaymath}

Let $A_i = A(e_i)$ be the open annulus in $\Omega(f)$ lying
over the edge $e_i$, and let $A$ be the annulus bounded
by $\bdry_+ A_n$ and $\bdry_+ A_0$.
Note that $f$ identifies the inner and outer boundaries
of $A$ bijectively, yielding a quotient torus
\begin{displaymath}
	T = \Abar/f .
\end{displaymath}
Since $f(w) = \lambda w$ in suitable local coordinates near $z$,
with $\lambda = f'(z)$, we have
\begin{displaymath}
	T \isom \cx/(2\pi i \zed \dirsum \log(\lambda) \zed).
\end{displaymath}
Let $B \subset T$ be the annulus homotopic to $A$ that is covered by
\begin{displaymath}
	\{w \st 0 < \Re(w) < \log |\lambda|\} \subset \cx.
\end{displaymath}
Since $\bdry B$ is geodesic, its modulus
\begin{displaymath}
	\mod(B) = \frac{\log |\lambda|}{2\pi} 
\end{displaymath}
is the maximum possible for any annulus for its homotopy class.

Applying Theorem \ref{thm:torus}, we obtain
\begin{displaymath}
	L(p,f) =
	\sum_1^n \delta(e_i) = 
	2\pi \sum_1^n \mod(A_i) 
	\le 2\pi \mod(B) = \log |f'(z)|
\end{displaymath}
as desired.
\qed

Let $O(1)$ denote a bound depending only on $\degf = \deg(f)$.

\begin{lemma}  \label{lemma:largeL}
If $L(p,F) \ge 6\pi \degf$, then $\log |f'(z)| \le L(p,F) + O(1)$.
\end{lemma}

\bold{Proof.} 
We continue the argument from the preceding proof.
Note that $\delta(e_i)$ is periodic, with period $n$,
for $i \leq n$.  Shifting indices, we can assume
$\delta(e_1) \ge \delta(e_i)$ for $1 < i \le n$
(and $\deg(e_n)=1$ as before).
Then the assumption $L(p,F) = \sum \delta(e_i) \ge 6\pi \degf$
implies $2\pi \mod(A_1) = \delta(e_1) \ge 6\pi$, and thus 
$\mod(A_1) \ge 3$ (using the fact that $n \le \degf$).
Thus we can apply the upper bound
of Theorem \ref{thm:torus} to obtain
\begin{displaymath}
	\log |f'(z)| = 2\pi \mod(B) 
	\le L(p,F) + 6\pi n (c_n/c_0)^2.
\end{displaymath}
Now recall that $f$ identifies the boundaries of $A$ and
expands metric $\rho = |\omega|$ by a factor of $\degf$.
Thus $(c_n/c_0) = \degf$, and therefore the defect
$6\pi n (c_n/c_0)^2$ is less than $6\pi \degf^3$, which depends only
on $\degf$.
\qed

\begin{lemma}  \label{lemma:smallL}  
If $L(p,F) < 6\pi \degf$, then $\log^+ |f'(z)| = O(1)$.
\end{lemma}

\bold{Proof.} 
Let $L_i = \sum_{i}^{i+n-1} \delta(e_i)$.
Then $L_i$ is monotone increasing, $\deg(e_i)L_i \leq L_{i+n} \le \degf L_i$, and
$L_i  < 6\pi \degf $ for $i \ll 0$. 
This implies we can find an index $j$ with $\deg(e_j) \ge 2$
and $1 \le L_j \le 6\pi \degf ^2 = O(1)$.
Now the monotone increasing sequence
\begin{displaymath}
	\deg(e_{j}), \deg(e_{j+n}), \deg(e_{j+2n}), \ldots
\end{displaymath}
can assume at most $\degf$ different values, so we can find a $k$ with $j\leq k \le j+\degf n$ such
that 
\begin{displaymath}
	2 \le \deg(e_k) = \deg(e_{k+n}) \le D.
\end{displaymath}
Since $L_{j+\degf n} \le \degf^\degf L_j$, we have $1 \le L_k \le O(1)$.

Now shift indices so that $k=0$; then $1 \le L_0 \le O(1)$.
Let $d = \deg(e_k)$.
Let $A_0,\ldots,A_{2n}$ be the annuli
lying over $e_0,\ldots,e_{2n}$.
Let $U_2 \subset U_1 \subset U_0$ be the disks in $\cx$ obtained by
filling in the bounded complementary components of $A_0$, $A_n$ and $A_{2n}$
respectively.  Then 
\begin{displaymath}
	f : U_1 \arrow U_0
\end{displaymath}
is a polynomial-like map of degree $d$,
and the fixed point $z$ of $f$ lies in $U_1$.
By construction, this polynomial-like map satisfies $U_2 = f^{-1}(U_1)$.
Since $\deg(e_0)=\deg(e_n)=d$, the critical points of
$f$ lie in $U_2$, and hence its critical values lie in $U_1$.

To control $\mod(U_0-\Ubar_1)$ and $\mod(U_0-\Ubar_2)$, we use
the flat metric $\rho = |\omega|$.
Note that $c_{2n}/c_0 \le \degf^2$,
since $f^2$ maps $\bdry_+ A_0$ onto $\bdry_+ A_{2n}$ and locally
expands the $\rho$-metric by a factor of $\degf^2$.
By the lower bound in Theorem \ref{thm:nest}, we have
\begin{displaymath}
	2\pi \mod(U_0-\Ubar_1) 
	\ge 2\pi \sum_{n+1}^{2n} \mod(A_i)
	= L_n \ge L_0 \ge 1,
\end{displaymath}
while the upper bound (together with $n = N(F) \le D-1$) yields:
\begin{displaymath}
	2\pi \mod(U_0-\Ubar_2) \le
	6\pi n (c_{2n}/c_0)^2 + 2\pi \sum_{1}^{2n} \mod(A_i) 
	\le 6\pi \degf^5 + L_0 + \degf L_0 = O(1) .
\end{displaymath}
Since the moduli of $U_0-\Ubar_1$ and $U_0-\Ubar_2$
are bounded above and below just in terms of $\degf$,
we have $\log^+ |f'(z)| = O(1)$ by Theorem \ref{thm:polylike}.
\qed

\bold{Proof of Theorem \ref{thm:mult}.}
Combine the results of Lemmas \ref{lemma:onesided}, \ref{lemma:largeL},
and \ref{lemma:smallL}.
\qed


\section{The moduli space of trees}
\label{sec:topology}

In this section we introduce the {\em geometric topology} on the 
moduli space $\cT_\degf$ of metrized polynomial-like trees of degree $\degf$.
Passing to the quotient projective space, we then show
$\proj\cT_\degf$ is compact and contractible (Theorem \ref{thm:compactcontractible}).

We also discuss the space $\cT_{\degf,1}$ of pointed trees and prove:

\begin{prop}  
\label{prop:lengthconv}
If $(T_n,\dist_n,F_n,p_n) \arrow (T,\dist,F,p)$ in $\cT_{\degf,1}$
and $F_n(p_n) = p_n$, then $F(p)=p$ and the 
translation lengths  satisfy
\begin{displaymath}
	L(p_n,F_n) \arrow  L(p,F).
\end{displaymath}
\end{prop}

\bold{The moduli space of trees.}
Let $\cT_\degf$ be the set of all equivalence classes of 
degree $\degf$ metrized polynomial-like trees $(T,\dist,F)$.
Trees $(T_1, \dist_1, F_1)$ and $(T_2,\dist_2, F_2)$ are equivalent 
if there exists an isometry $i: T_1\to T_2$ such that $i\circ F_1 = F_2\circ i$. 

There is a natural action of $\reals_+$ on $\cT_\degf$ which simply rescales
the metric $\dist$; the quotient projective space will be denoted $\PT_\degf$.
A tree is {\em normalized} if $\dist(v_0,J(F)) = 1$, where $v_0$ is the base
of the tree.
The normalized trees form a cross-section to the projection
$\cT_\degf \arrow \PT_\degf$.

\bold{Strong convergence.}
Let $v_i \mem V(T)$, denote the unique vertex at combinatorial
height $h(v_i)=i \ge 0$, and let $T(k) \subset T$ denote the finite subtree
spanned by the vertices with combinatorial height $-kN(F) \le h(v) \le kN(F)$.
Recall that $N(F)$ is the number of disjoint grand orbits of vertices, as introduced
in Lemma \ref{lemma:h}.

We say a sequence $(T_n,\dist_n,F_n)$ in $\cT_\degf$
{\em converges strongly} if:
\begin{enumerate}
	\item
The distances $\dist_n(v_0,v_i)$ converge for $i=1,2,\ldots,\degf$;
	\item
We have $\lim \dist_n(v_0,v_\degf) > 0$; and
	\item
For any $k>0$ and $n > n(k)$, there is a simplicial isomorphism
$T_n(k) \isom T_{n+1}(k)$ respecting the dynamics.
\end{enumerate}
The last condition implies $N(F_n)$ is eventually constant.

\begin{lemma}
\label{lemma:compactness}
Any sequence of normalized trees in $\cT_\degf$ has a
strongly convergent subsequence.
\end{lemma}

\bold{Proof.} 
In a sequence of normalized trees, $\dist_n(v_0,v_i) \le \degf^i$ and
$\dist_n(v_0,v_\degf) \ge 1$, so the first two properties of strong
convergence hold along a subsequence.
The number of vertices in $T_n(k)$
is bounded in terms of $\degf$ and $k$, so the third property
holds along a further subsequence.
\qed

\bold{Limits.}
Suppose $(T_n,\dist_n,F_n)$ converges strongly.  Then there is a unique pointed simplicial complex $(T',v_0)$ with dynamics $F' : T' \arrow T'$ such that $T_n(k) \isom T'(k)$ for all $n > n(k)$, and the simplicial isomorphism respects the dynamics.  It is possible, however, that certain edge lengths of $T_n$ tend to 0 in the limit; this happens when the grand orbits of critical points collide.
Our assumptions therefore yield only a pseudo-metric $\dist'$ on $T'$
as a limit of the metrics $\dist_n$.
Let $(T,\dist,F)$ be the metrized dynamical system obtained by
collapsing the edges of length zero to points.

\begin{lemma}   \label{lemma:degreelimit}
Suppose $(T_n,\dist_n,F_n)$ converges strongly.
The limiting triple $(T,\dist,F)$ is a metrized polynomial-like tree.
\end{lemma}

\bold{Proof.} 
Let the vertices $V(T)$ be the grand orbits of its branch points.
Since $\lim \dist_n(v_0,v_\degf) > 0$, $V(T)$ is nonempty, and it is easy to see that
$T$ has the structure of a locally finite simplicial tree, and $F : T \arrow T$ is
a branched cover.
We must show $T$ has a compatible degree function.

To define this, pass to a subsequence such that
for each $k$ the degree function of $T_n$ restricted to $T_n(k)$ stabilizes as $n \arrow \infty$.
This defines a degree function $\deg' : E(T') \arrow \natls$ on the simplicial limit $T'$
compatible with $F'$.  

Note that $T'$ may have vertices of valence two whose grand orbits under $F'$ contain
no branch points.  These vertices arise when the critical point that used to label them
no longer escapes.  Since they have valence two, the value of $\deg'$ is the same on
both their adjacent edges.  We can thus modify the simplicial structure on $T'$
by removing all such vertices, and maintain a compatible degree function by taking
its common value on all edges that are coalesced.

With this modified simplicial structure on $T'$,
the natural collapsing map $T' \arrow T$ is simplicial.
We define $\deg : E(T) \arrow \natls$ by $\deg(e)=\deg'(e')$ for
the unique edge $e'$ lying over $e$, and for $v \mem V(T)$ define
$\deg(v)=\deg(e)$ where $e$ is the upper edge of $v$.
It is then straightforward to check that the resulting degree function
is compatible with $F : T \arrow T$.
\qed

\bold{The geometric topology.}
The {\em geometric topology} on $\cT_\degf$ is the unique 
metrizable topology satisfying
\begin{displaymath}
	(T_n,\dist_n,F_n) \arrow (T,\dist,F)
\end{displaymath}
whenever $(T_n,\dist_n,F_n)$ is strongly convergent and $(T,\dist,F)$ is defined as above.
In \S\ref{sec:GH}, we show that the geometric topology coincides with the Gromov-Hausdorff topology on {\em pointed dynamical metric spaces}; in particular, we describe there a basis of open sets for the topology.  

Lemma \ref{lemma:compactness} immediately implies:

\begin{theorem}  
\label{thm:projcompact}
The space $\PT_\degf$ is compact in the quotient geometric topology.
\end{theorem}

\bold{Iteration.}  For each $(T,\dist, F)\in\cT_\degf$, its $n$-th iterate $(T,\dist, F^n)$
is a metrized polynomial-like tree of degree $\degf^n$.  
Define
	$$i_n: \cT_\degf \to \cT_{\degf^n}$$
by $(T,\dist, F)\mapsto (T, \dist, F^n)$.  It is useful to observe:

\begin{lemma}  \label{lemma:iteration}
The iterate maps $i_n$ are continuous in the geometric topology.
\end{lemma}

\bold{Proof.}
It suffices to consider sequences $(T_m,\dist_m, F_m)$ converging strongly
to $(T,\dist, F)$ in $\cT_\degf$.  For each $k>0$, any simplicial isomorphism 
$s: T_m(k) \to T'(k)$ such that $s\compos F_m = F'\compos s$
will also satisfy $s\compos F_m^n = (F')^n\compos s$.  Therefore, the 
sequence $(T_m, \dist_m, F_m^n)$ converges strongly to $(T, \dist, F^n)$.
\qed

Next we establish:

\begin{theorem}  \label{thm:contractible}
The space $\proj\cT_\degf$ is contractible.
\end{theorem}

The proof is based on a natural construction which {\em accelerates} the rate of escape of critical points in a tree.  A version of the following result appears as Theorem 7.5 in \cite{Emerson:trees}.

\begin{theorem}
\label{thm:subtree}
Let $(T,\dist,F)$ be a metrized polynomial-like tree, and let $S \subset T$ be a
forward-invariant subtree. 
Then $F|S$ can be extended to a unique metrized polynomial-like tree $(T',\dist',F')$
with the same degree function on $S$, and whose critical points all lie in $\Sbar$.
\end{theorem}

We emphasize that the subtree $S$ can have endpoints, and that these endpoints
need not coincide with vertices of $T$.  The degree of a terminal edge of $S$ 
is defined to be the degree of the edge of $T$ which contains it.

\bold{Proof.} 
The characterization of critical points in Lemma \ref{criticalpoints} requires
that all edges in $T'\setminus S$ have degree 1.
The tree $T'$ and the map $F':T'\arrow T'$ will be defined inductively
on (descending) height, uniquely determined
by the conditions that each added edge has degree 1 and that
(\ref{eq:degineq}) and (\ref{eq:degeq})
are satisfied at all vertices of $T'$.

Let $p$ be a point of maximal height in $\closure{T\setminus S}$; set 
$T' = \closure{S}$, $F'|T' = F|S$, and $\dist'|T' = \dist|S$. 
Then $p$ is a highest point in $T'$ such that either 
(a) $F'(p)$ lies in the interior of an edge of $T'$, or (b) $F'(p)$ is a vertex and the 
local degree condition (\ref{eq:degeq}) for $F'|T'$ is not satisfied at $p$. 

In case (a), the point $p$ belonged to the interior of an edge $e$ of $T$.  
We make $p$ into a vertex of degree $\deg(e)$.  
Extend $(T',\dist',F')$  below $p$ down to height $H(p)/d$
to be a local homeomorphism, defining $\dist'$ so that 
$\dist'(e') = \dist'(F(e'))/d$ on each new edge $e'$.  Assigning new edges 
degree 1, the conditions
(\ref{eq:degineq}) and (\ref{eq:degeq}) will both be satisfied at $p$.  
Note that the degree conditions are
always satisfied at vertices where $F'$ is a local homeomorphism 
and all adjacent edges have degree 1.

In case (b), define $(T',\dist',F')$ in a neighborhood of $p$ by adding enough 
new edges of degree 1 below $p$ so that 
the local degree condition (\ref{eq:degeq}) is satisfied with degree $\deg(p)$.    
Again, we can define $(T',\dist',F')$ on the 
added edges and vertices of $T'$ below $p$ 
down to height $H(p)/d$ so that $F'$ is a local homeomorphism and 
$\dist'(e') = \dist'(F(e'))/d$ on all new edges $e'$.  
Condition (\ref{eq:degineq}) will be automatically satisfied at $p$ because 
it is satisfied at $p$ for $(T,F)$ and 
the right-hand side can only decrease with the replaced edges of degree 1.  

There are only finitely many endpoints or vertices $x\in T'$ with height
$H(p)/d < H(x) \leq H(p)$ where (a) or (b) is satisfied, 
and we repeat the above construction for each of 
these points.  We then may proceed by induction on height 
of vertices where the local degree is not well-defined, until we have completed
the construction of $(T',\dist', F')$.  
\qed

\bold{Escaping trees.}  
A metrized polynomial-like tree $(T,\dist, F)$ is {\em escaping} if there 
are no critical points in $J(F)$. 

\begin{cor}  
\label{cor:denseescape}
Escaping trees are dense in the spaces $\cT_\degf$ and $\proj\cT_\degf$.
\end{cor}

\bold{Proof.}  
Let $(T,\dist,F)$ be a metrized polynomial-like tree with height function $H:\Tbar\to [0,\infty)$.  
For each $\epsilon>0$, let $S_\epsilon\subset T$
be the subtree of all points with height $\geq \epsilon$.  By Theorem \ref{thm:subtree},
we can extend $F|S_\epsilon$ uniquely so that all critical points are contained in $S_\epsilon$
to obtain $(T_\epsilon, \dist_\epsilon, F_\epsilon)$.  Letting $\epsilon\to 0$, we have 
 $$(T_\epsilon, \dist_\epsilon, F_\epsilon) \to (T,\dist,F)$$
in the geometric topology.
\qed

\bold{Proof of Theorem \ref{thm:contractible}.}  
Identify $\proj\cT_\degf$ with the subset of normalized trees in $\cT_\degf$.  
For each normalized tree $(T, \dist, F)$ with height function 
$H: \Tbar\to [0,\infty)$ and each $t\in [0,1]$, consider
the forward-invariant subtree
 $$S_t = \{x\in \Tbar: H(x) \geq t\}.$$
By Theorem \ref{thm:subtree},
there is a unique metrized polynomial-like tree $(T_t, \dist_t, F_t)$ extending 
$F|S_t$ and the local degree function on $S_t$ so that all critical points 
belong to $S_t$.    

Define
	$$R:\proj\cT_\degf\times [0,1]\to \proj\cT_\degf$$
by $((T, \dist,  F), t) \mapsto (T_t, \dist_t, F_t)$.   Then $R(\,\cdot\, , 0)$ is 
the identity, and $R(\,\cdot\, , 1)$ is the constant map sending all
trees to the unique normalized tree $(T_1, \dist_1, F_1)$ 
with all critical points at the base $v_0$.   Note that $R((T_1, \dist_1, F_1),t) =
(T_1, \dist_1, F_1)$ for all $t$.  
It remains to show that $R$ is continuous.

Fix $(T, \dist, F)$, $t\in [0,1]$, a sequence $(T_n, \dist_n, F_n)$
of normalized trees converging strongly to $(T,\dist,F)$, and 
a sequence $t_n\to t$.    Because the number of critical points 
(and thus their grand orbits) is finite,
we may pass to a subsequence so that the subtrees
$S_{t_n}\subset T_n$ are simplicially isomorphic (respecting dynamics)
for all $n>>0$.  
The isomorphisms can be extended to $T_{n,t_n}\isom T_{n+1, t_{n+1}}$ 
using the construction of $F_{n,t_n}$ as a local homeomorphism below 
$S_{t_n}$.  Therefore, the image sequence $R((T_n, F_n), t_n)$
converges strongly.  The limit clearly coincides with $T$ above height 
$t$.  Because the degree functions converge, it must have all edges of 
degree 1 below height $t$.  By the uniqueness in Theorem \ref{thm:subtree},
the limit must be $(T_t, \dist_t, F_t)$.
\qed

\bold{Proof of Theorem \ref{thm:compactcontractible}.}
Combine Theorems \ref{thm:projcompact} and \ref{thm:contractible}.
\qed

\bold{Pointed trees.} 
A {\em pointed tree} is a quadruple $(T,\dist, F,p)$  
where $(T,\dist,F) \mem \cT_\degf$ and $p \mem \Tbar$.
Let $\cT_{\degf,1}$ denote the set of isometry classes pointed trees of degree $\degf$.

Let $p(k) \mem T(k)$ denote the image of $p \mem T$ under the nearest-point
retraction $\Tbar \arrow T(k)$.
We say a sequence $(T_n, \dist_n,F_n, p_n)$ in $\cT_{\degf,1}$ {\em converges strongly} if 
\begin{enumerate}
\item	$(T_n, \dist_n,F_n)$ converges strongly; 
\item 	$\dist_n(v_0, p_n)$ converges to a finite limit as $n \arrow \infty$; and
\item	for all $k>0$ and all $n>n(k)$, there exists a simplicial isomorphism of pointed spaces
		$(T_n(k),p_n(k))\isom (T_{n+1}(k),p_{n+1}(k))$ respecting the dynamics.
\end{enumerate}
In this case the pointed isomorphisms on finite trees determine a natural
pointed limit $(T,\dist,F,p)$,
and we define the {\em geometric topology} on $\cT_{\degf,1}$ by requiring that
$(T_n,\dist_n,F_n,p_n) \arrow (T,\dist,F,p)$ for every strongly convergent sequence.
(Similar definitions can be given for $\cT_{d,m}$, $m>1$.)

\bold{Continuity of translation lengths.}
This space of pointed trees is useful for tracking periodic points and critical points.
For example, it is straightforward to verify:

\begin{prop}
\label{prop:cps}
The set of normalized pointed trees $(T,\dist,F,p)$ such that $p$ is a critical point of $F$
is compact in $\cT_{\degf,1}$.
\end{prop}
We can now establish continuity of translation lengths.

\bold{Proof of Proposition \ref{prop:lengthconv}.}
It is enough to treat the case where 
\begin{displaymath}
	(T_n,\dist_n,F_n, p_n) \arrow (T,\dist,F,p)
\end{displaymath}
strongly; then clearly $F(p)=p$.
If $p$ is not a critical point of $F$, then there is a $k>0$ such that
$T$ has no critical points below $p(k)$.
By Proposition \ref{prop:cps}, $F_n$ has no critical points below $p_n(k)$ for $n \gg 0$,
and thus
\begin{displaymath}
	L(F_n,p_n) = \delta_n(p_n(k),F_n(p_n(k))) .
\end{displaymath}
By geometric convergence, the metric $\dist_n|T_n(k)$ converges to $\dist|T(k)$, and similarly
for the degree function; thus the corresponding modulus metrics also satisfy $\delta_n \arrow \delta$ on 
finite subtrees, and hence
\begin{displaymath}
	\delta_n(p_n(k),F_n(p_n(k))) \arrow \delta(p(k),F(p(k)))  = L(F,p).
\end{displaymath}

On the other hand, if $p$ is a critical point then $L(F,p)=0$ and hence
$\delta(p(k),F(p(k))) \arrow 0$ as $k \arrow \infty$.
By geometric convergence, $p_n(k)$ is also moved a small amount by
$F_n$ when $n \gg 0$, and thus $L(p_n,F_n) \arrow 0$.
\qed


\section{Continuity of the quotient tree}
\label{sec:continuity}

In this section we study the map from the moduli space of polynomials
to the moduli space of trees, and establish:

\begin{theorem}
\label{thm:cpe}
The map $\tau: \MPoly_\degf^* \arrow \cT_\degf$ is continuous, proper, 
and equivariant with respect to the action of $\reals_+$ by stretching of
polynomials and by metric rescaling of trees.
\end{theorem}

\noindent
This gives Theorem \ref{thm:quotmap} apart from surjectivity, which will be
established in \S\ref{sec:realization}.

\bold{The moduli space of polynomials.}
Let $\MPoly_\degf = \Poly_\degf/\Aut(\cx)$ be the moduli space
of polynomials of degree $\degf\geq 2$.  Every polynomial is conjugate to one which
is monic and centered, {\em i.e.} of the form 
$$f(z) = z^\degf + a_{\degf-2} z^{\degf-2} + \cdots a_1 z + a_0$$
with coefficients $a_i\in\cx$, 
and thus $\MPoly_\degf$ is a complex orbifold finitely covered by $\cx^{\degf-1}$.

The escape-rate function of a polynomial satisfies   $G_{AfA^{-1}}(Az) = G_f(z)$
for any $A \mem \Aut(\cx)$. Consequently, the
{\em maximal escape rate}
\begin{displaymath}
	M(f)=\max\{G_f(c): f'(c)=0\}
\end{displaymath}
is well-defined on $\MPoly_\degf$.  
The open subspace $\MPoly_\degf^*$ where $J(f)$ is disconnected coincides with
the locus $M(f)>0$.

By Branner and Hubbard \cite[Prop 1.2, Cor 1.3, Prop 3.6]{Branner:Hubbard:cubicsI}
we have:

\begin{prop}  
\label{prop:BHcontinuous}
The escape-rate function $G_f(z)$ is continuous in both $f\in \Poly_\degf$ 
and $z\in\cx$.  
\end{prop}

\begin{prop}  
\label{prop:BHproper}
The maximal escape rate  $M: \MPoly_\degf^* \arrow (0,\infty)$ is 
proper and continuous.
\end{prop}

\bold{Stretching.}
The {\em stretching deformation} associates to any polynomial $f(z)$ of degree $\degf>1$
a 1-parameter family of topologically conjugate polynomials $f_t(z)$, $t \mem \reals_+$.
To define this family, note that the Beltrami differential defined by
\begin{displaymath}
	\mu = \frac{\omegabar}{\omega}
\end{displaymath}
on the basin of infinity, where $\omega = 2 \del G_f$,
and $\mu=0$ elsewhere, is invariant under $f$.
Consequently, if we let $\phi_t : \cx \arrow \cx$ be a smooth family of quasiconformal maps solving the
Beltrami equation
\begin{displaymath}
	\frac{d\phi_t/d\zbar}{d\phi_t/dz} = 
	\frac{t-1}{t+1} \mu,
\end{displaymath}
$t \mem \reals_+$, then 
\begin{displaymath}
	f_t = \phi_t \compos f \compos \phi_t^{-1}
\end{displaymath}
is a smooth family of polynomials with $f_1 = f$.
The maps $\phi_t(z)$ behave like $(r,\theta) \mapsto (r^t,\theta)$ near infinity, and thus
the corresponding Green's functions satisfy
\begin{displaymath}
	G_{f_t}(\phi_t(z)) = t G_f(z) 
\end{displaymath}
(compare \cite[\S 8]{Branner:Hubbard:cubicsI}).  
Together with Proposition \ref{prop:BHproper}, this implies:  

\begin{prop}
\label{prop:stretchpoly}
For any polynomial $f$ with disconnected Julia set, the stretched polynomials $f_t$
determine a smooth and proper map $(0,\infty)\to \MPoly_\degf^*$.  
\end{prop}

In addition:  

\begin{prop}
\label{prop:stretch}
The quotient tree for the stretched polynomial $f_t$ is obtained from the
quotient tree $(T,\dist, F)$ for $f$ by replacing the height metric $\dist(x,y)$ with $t\dist(x,y)$.
\end{prop}
Note that there is also a twisting deformation, using $i \mu$, which does not change the quotient tree for $f$.

\bold{Proof of Theorem \ref{thm:cpe}.}
Equivariance of $\tau : \MPoly_\degf^* \arrow \cT_\degf$ 
with respect to stretching is Proposition \ref{prop:stretch}.

To prove continuity,
suppose $[f_n] \arrow [f]$ in $\MPoly_\degf^*$.
Lift to a convergent sequence $f_n \arrow f$ in $\Poly_\degf$.
Since $M(f_n) \arrow M(f)>0$
we can pass to a subsequence so the corresponding trees 
$(T_n,\dist_n, F_n)$ converge strongly to $(T,\dist,F) \mem \cT_\degf$.
It suffices to show that $(T,\dist,F)$ is isometric to the tree for $f$.

By the definition of strong convergence, we have a limiting simplicial
tree map $F' : T' \arrow T'$ with a pseudo-metric $\dist'$,
and simplicial isomorphisms $T'(k) \isom T_n(k)$ for all $n>n(k)$, respecting the dynamics (see \S\ref{sec:topology} where the geometric topology is introduced).
Fix $k>0$, and recall that the Green's function $G_n$ for $f_n$ factors through $T_n$.
Moreover the subtree $T_n(k)$ corresponds to the compact region
\begin{displaymath}
	\Omega_k(f_n) =
	\{z \mem \cx \st \degf^{-k}M(f_n) \le G(z) \le \degf^k M(f_n) \}.
\end{displaymath}
Thus the vertices of $T'(k)$ label components of the critical level sets of $G_n$ in this range for all $n$ sufficiently large.  By Proposition \ref{prop:BHcontinuous},
$G_n$ converges uniformly on compact sets to the Green's function $G$ for $f$.
Thus $\Omega_k(f_n)$ converges to $\Omega_k(f)$, and
we obtain a corresponding labeling of the critical level sets of $G$ by $T'(k)$ (though multiple vertices can label the same component of a level set).
The distance $\dist'(v_1,v_2)$ between consecutive
vertices in $T'$ encoding level sets $L_1$ and $L_2$ is given simply by
$|G(L_1)-G(L_2)|$.  It follows that 
$(T,\dist,F)$ is exactly the quotient tree for $f$, and thus $\tau$ is continuous.

Finally Proposition \ref{prop:BHproper} implies that
$\tau$ is proper, since
$M(f) = \dist(v_0,J(F))$ is bounded above and below on any compact subset of $\cT_\degf$.
\qed

\bold{Remark:  planar embeddings.}
Topologically, the level sets of the Green's function of $f(z)$ 
are graphs embedded in $\cx$. 
These planar graphs are not always uniquely determined by the tree of $f$,
and thus
the map $\tau : \MPoly_\degf^* \arrow \cT_\degf$ can have disconnected fibers.
In the simplest examples, different graphs correspond to different choices
for a primitive $n$th root of unity, suggesting a connection with Galois theory and
{\em dessins d'enfants}; cf.  \cite{Pilgrim:dessins}


\section{Polynomials from trees}
\label{sec:realization}

In this section we prove:
\begin{theorem}  \label{thm:realization}
Any metrized polynomial-like tree $(T,\dist,F)\in\cT_\degf$ can be realized by
a polynomial $f$.
\end{theorem}

Together with Theorem \ref{thm:cpe},
this completes the proofs of Theorems \ref{thm:Tsurjective}
and \ref{thm:quotmap} of the Introduction.

\bold{Permutations.}
A partition $P$ of $D \ge 1$ is an unordered sequence of 
positive integers  $(a_1,\ldots,a_m)$ such that
$D = a_1+ \cdots + a_m$.
A partition $P$ of $D$ determines a conjugacy class $S_D(P)$ in 
the symmetric group $S_D$, consisting of 
all permutations which are products of $m$ disjoint cycles with lengths
$(a_1, a_2, \ldots, a_m)$.

Let $c(P) = D-m = \sum (a_i-1)$.  In our application to branched coverings,
$c(P)$ will count the number of critical points coming from the blocks of $P$.

\begin{prop} 
\label{prop:hurwitz} 
Let $P_1,\ldots,P_n$ be partitions of $D$ such that $\sum_1^n c(P_i) = D-1$.
Then there exist permutations $\sigma_1,\ldots,\sigma_n$ in the corresponding
conjugacy classes of $S_D$, such that
$\sigma_1 \cdots \sigma_n = (123\ldots D)$.
\end{prop}

\bold{Proof.} 
First note that if $P = (a_1,\ldots,a_m)$ and $c(P) < D/2$, 
then $m > D/2$ and thus $a_i=1$ for some $i$.

We proceed by induction on $D$, the case $D=1$ being trivial.
Assume the result for $D' = D-1$.
Let us order the partitions $P_i$ and their entries $(a_1, \ldots, a_m)$ 
so that $c(P_1) \ge c(P_i)$ and $a_1 \geq a_i$ for all $i$.  
Then $c(P_1)>0$ so $a_1>1$, 
and $c(P_i) < D/2$ for $i>1$, so each of these partitions has at
least one block of size $1$.

Let $P_1' = (a_1-1,a_2,\ldots a_m)$,
and define $P_i'$, $i>1$ by discarding a block of size $1$ from $P_i$.
Then $P_1',\ldots,P_n'$ are partitions of $d'$ satisfying
$\sum c(P_i') = D'-1 = D-2$.  By induction there are permutations
$\sigma_i' \mem S_{D-1}$
corresponding to $P_i'$ whose product is the cycle $(123\ldots D')$.

We can assume that $1$ belongs to the cycle of length $(a_1-1)$ for
$\sigma_1'$.  Then $\sigma_1 = (1D)\sigma_1' \mem S_D$ has a cycle of length
$a_1$ and overall cycle structure given by $P_1$.
Taking $\sigma_i = \sigma_i'$ for $i>1$ (under the natural inclusion
$S_{D-1} \includesin S_D$), we find $\sigma_i$ has an additional cycle of
length $1$ and hence it lies in the conjugacy class $S_D(P_i)$.
Finally we have
\begin{displaymath}
	\sigma_1 \cdots \sigma_n = (1D)\sigma_1' \cdots \sigma_n' =
		(1D)(123\ldots (D-1)) = (123\ldots D).
\end{displaymath}
\qed

\bold{Branched coverings.}
Suppose $f: X\arrow Y$ is a degree $D$ branched covering of Riemann surfaces.  
Given $y \mem Y$, the {\em branching partition} of $f$ over $y$ is
the partition of $D$ given in terms of the local degree
of $f$ at each of the preimages $f^{-1}(y) = (x_1,\ldots,x_m)$ by
\begin{displaymath}
	P(f,y) = (\deg(f,x_1),\ldots,\deg(f,x_m)).
\end{displaymath}
The quantity $c(P(f,y))=D-m$ is the number of critical points in the fiber $f^{-1}(y)$,
counted with multiplicities.

Suppose now that $f:\cx\to\cx$ is a polynomial of degree $\degf$, with critical values $\{p_1, \dots, p_n\}$.  Choose a basepoint $b$ which is not a critical value of $f$.  Then the fundamental group $\pi_1(\cx\setminus\{p_1, \ldots, p_n\}, b)$ acts by permutations on the fiber $f^{-1}(b)$.  If $\sigma_i$ denotes the permutation induced by a loop around $p_i$, then up to relabeling, the product $\sigma_1\sigma_2\cdots\sigma_n$ is equal to the permutation $(123\cdots D)$ which is the permutation induced by a loop around $\infty$.  

Let $(T,\dist, F)$ be a polynomial-like tree, and let $v$ be a vertex of $T$.
A polynomial $f:\cx\arrow \cx$ of degree $\deg(v)$ has 
the {\em branching behavior of $(T,F,v)$ over $p_1, \ldots, p_n\in\cx$} if
there is an ordering of the 
lower edges $e_1, \ldots, e_n$ of $T$ at  $F(v)$ such that
\begin{displaymath}
	P(f,p_i) = (\deg(e) : e\in E_v,  F(e)=e_i)
\end{displaymath}
for $i=1,\ldots,n$.
If the critical multiplicity 
$$m(v) = 2\deg(v)-2 - \sum_{e\in E_v} (\deg(e)-1)$$
is non-zero, then  $f$ will have critical values outside the set $\{p_1, \ldots, p_n\}$.  

\begin{prop} \label{localexist}
Let $(T,\dist,F)$ be a metrized polynomial-like tree, 
$v$ a vertex of $T$, and $n$ the number of lower edges at $F(v)$.   
For any set of distinct points $\{p_1, \ldots, p_n, q\}$ in $\cx$, 
there exists a polynomial
of degree $\deg(v)$ with the branching behavior of $(T,F,v)$ over 
$p_1, \ldots, p_n$ and all critical values contained in 
$\{p_1, \ldots, p_n, q\}$.   
\end{prop}

\bold{Proof.} 
Let $e_1, \ldots, e_n$ be the lower edges of $T$ at $F(v)$.
For each $e_i$, its set of preimages in $E_v$ determines
the partition $P_i$ of $\deg(v)$ given by 
$(\deg(e): F(e)=e_i)$.  Let $Q$ be the partition
$(m(v) + 1, 1, \ldots ,1)$ of $\deg(v)$.  Then $c(Q) + \sum c(P_i) = \deg(v)-1$.

By Proposition \ref{prop:hurwitz}, there 
exist permutations $\sigma_1, \ldots, \sigma_n, \sigma_q$ in 
the corresponding conjugacy 
classes of the symmetric group $S_{\deg(v)}$ 
with product $\sigma_1 \cdots \sigma_n\sigma_q = (12\ldots \deg(v))$.
The representation 
	$$\pi_1(\cx \setminus \{p_1, \ldots, p_n, q\}) \arrow S_{\deg(v)}$$
which associates to each generating loop the permutation $\sigma_i$ 
or $\sigma_q$
determines a holomorphic branched
covering $f: \chat\to\chat$, with branching partitions $P(f,p_i) = P_i$,
$P(f,q) = Q$ and $P(f,\infty) = (\deg(v))$.  In particular, $f$ is
totally ramified over $\infty$.  Choosing coordinates on the domain so that
$f(\infty)=\infty$, we find that $f$ is a polynomial with the required branching behavior.
\qed

\bold{Proof of Theorem \ref{thm:realization}.}   
We will first prove the realization theorem in the {\em escaping} case, where
$(T,\dist,F)$ has no critical points in its Julia set $J(F)$.
The general case will follow by density of escaping trees and a compactness argument, using the continuity of $\tau: \MPoly_\degf^*\to\cT_\degf$ (Theorem \ref{thm:cpe}).  

Let $(T,\dist,F)$ be an escaping tree of degree $\degf$.   For each vertex $v$ of $T$, we will use Proposition \ref{localexist} to construct a local polynomial realization
  $$f_v:\cx_v \arrow \cx_{F(v)},$$ 
together with a foliation $\cF_v$ of $\cx_v$.  The foliation will have the following structure:  its leaves are the level sets of a subharmonic function $G_v:\cx_v\to [-\infty, \infty)$ with $\Delta G_v = \sum c_i \delta_{\zeta_i}$ for a finite collection of points $\zeta_i$ in bijective correspondence with the lower edges of $v$, the level set $L_v = \{G_v = 0\}$ is connected, and the connected components of $\{G_v < 0\}$ are topological disks each containing a unique $\zeta_i$.  We require the compatibility condition 
\begin{equation} \label{eq:Gcomp}
	G_v(z) = G_{F(v)} (f_v(z))/ \deg(v),
\end{equation}
so that $f_v$ pulls back the foliation $\cF_{F(v)}$ to the foliation $\cF_v$, taking the {\em central leaf} $L_{F(v)}$ to the central leaf $L_v$.  We then glue the local realizations along leaves of the foliations to obtain a polynomial $f$ such that $\tau(f) = (T,\dist,F)$.

\bold{The local models.}
Fix a vertex $v$ with critical multiplicity $m(v)$, and assume that $F(v)$ 
is a vertex of valence 2;  
this is always the case if the combinatorial height of $v$ is $\geq 0$.   
Mark the point $p_1 = 0$  in $\cx_{F(v)}$.  Let $G_{F(v)}(z) = \log|z|$; 
the associated foliation of $\cx_{F(v)}$ is by circles $|z|=c$ with the unit circle as central leaf.  For $m(v)\not=0$, let
$q$ be a point on the unit circle.
Let $f_v:\cx_v\to\cx_{F(v)}$ be any polynomial guaranteed
by Proposition \ref{localexist} with the branching behavior of $(T,F,v)$ 
over $p_1$ and critical values $\{p_1,q\}$.  Define $G_v$ on $\cx_v$ by the compatibility condition (\ref{eq:Gcomp}).  
The foliation by circles $|z|=c$ in $\cx_{F(v)}$
pulls back to a (singular) foliation of $\cx_v$: the preimages of the circle $|z|=c$
with $c\not=0,1$ are topological circles, and the central leaf is a connected degree $\deg(v)$ branched cover of the unit circle, branched over one point with multiplicity $m(v)$.   The preimages of the marked point $p_1$ are indexed by the edges below $v$.
For the case $m(v)=0$, we can take $f_v(z) = z^{\deg(v)}$.  

We complete the definitions of the local realizations by induction.  Assume
that $f_v:\cx_v \arrow \cx_{F(v)}$ has been defined and the foliation with distinguished 
central leaf has been specified on the domain.  There is also a marked set
of points in $\cx_v$ corresponding to the lower edges adjacent to $v$.  For each vertex $v'$ such that 
$F(v')=v$, we use Proposition \ref{localexist} to 
define the polynomial $f_{v'}$ with the branching behavior of $(T,F,v')$ over the
the marked points 
in $\cx_v$ with branch point of multiplicity $m(v')$ over an arbitrary point $q$ on the central 
leaf.

\bold{Cutting and pasting.}
For each vertex $v$, we define a Riemann surface with boundary $S_v\subset \cx_v$
according to the data of $(T,\dist,F)$.
For vertices connected by an edge, we will glue the associated surfaces so that
the local maps match up.  

Let $v_0$ be the base of $T$.  Consider the consecutive
vertices $v_0,  v_1, \ldots, v_n = F(v_0), v_{n+1} = F(v_1)$, 
bounding edges $e_0, e_1, \ldots, e_n$ of lengths
$l_0, \ldots, l_n$, where $l_n=\degf l_0$, in the height metric $\dist$.   
For each $i=1, \dots, n$, let  
 $$S_{v_i} = \{e^{-l_{i-1}} \leq |z| \leq e^{l_i}\} \subset \cx_{v_i}$$
with central leaf $\{|z|=1\}$.   
For each $i$, we identify the outer boundary of $S_{v_i}$ with the inner boundary
of $S_{v_{i+1}}$ via an isometry with respect to 
the metric $|dz/z|$ to form a cylinder; the twist parameters are free.
Because the leaves $\{|z|=c\}$ are extremal curves of these annuli, 
the central leaves of $S_{v_i}$ and $S_{v_{i+1}}$ 
bound an annulus of modulus exactly $(l_i/4\pi) + (l_i/4\pi) = l_i/2\pi$.  
For the vertex $v_0$, let $S_{v_0} = f_{v_0}^{-1}(S_{v_n})\subset \cx_{v_0}$, 
and glue the outer boundary of $S_{v_0}$ to the inner boundary of $S_{v_1}$.  
By construction, the modulus of the annulus bounded by the central leaves of 
$S_{v_0}$ and $S_{v_1}$ is therefore $l_n/(4\pi \degf) + l_0/4\pi = l_0/2\pi$.   
The holomorphic functions $f_{v_i}$ and $f_{v_{i+1}}$ extend
across the common boundary of $S_{v_i}$ and $S_{v_{i+1}}$ for all 
$i = 0, \ldots,n$.

We are now set up for an inductive construction.  Suppose that $v$ and $w$ 
are
two vertices connected by an edge, and suppose we have defined $S_v$,  
$S_w$, and the gluing between them.   
Let $v'$ and $w'$ be adjacent vertices such that $F(v')=v$ and $F(w')=w$.  
Set $S_{v'} = 
f_{v'}^{-1}(S_v)$ and $S_{w'} = f_{w'}^{-1}(S_w)$.  Let $e$ be the edge
connecting $v'$ and $w'$.  There are exactly $\deg(e)$ ways to glue 
$S_{v'}$ and $S_{w'}$ so that the maps $f_{v'}$ and $f_{w'}$ extend 
across the common boundary; make any of these choices.  

It remains to consider the edges of combinatorial height $> N(F)$.  Suppose $v$ 
and $w$ are
vertices connected by an edge $e$ of degree $\degf$, and let $V$ and $W$ be their 
images under $F$.  Let $S_V = f_v(S_v)\subset \cx_V$ and $S_W = 
f_w(S_w)\subset \cx_W$.  In this setting, there is a unique gluing of 
$S_V$ and $S_W$ so that the maps $f_v$ and $f_w$ extend continuously across the
common boundary of $S_v$ and $S_w$.  

\bold{The result of the inductive construction.}
We have produced a holomorphic map $f: S\arrow S$ 
on a planar Riemann surface $S$ equipped with a foliation such that 
$F: T\arrow T$ is the quotient of $f:S\arrow S$ by this foliation.  Furthermore, 
to every edge $e$ in $T$ is associated an annulus $A_e\subset S$ with modulus 
satisfying $\mod (A_e) = \mod(f(A_e))/\deg(e)$.  If $e$ is an edge 
contained in the path $[v_0,\infty)$, 
then $\dist(e) = 2\pi \mod(A_e)$.  

\bold{The map $f$ extends to a polynomial.}
Since $S$ is planar, there exists a holomorphic embedding $S \includesin \cx$
sending the unique isolated end of $S$ to infinity \cite[\S 9-1]{Springer:book:RS}.
Because $(T,\dist,F)$ is an escaping metrized polynomial-like tree, there is a height $\epsilon>0$ so that all edges of height $<\epsilon$ have degree 1.  These edges
give chains of disjoint annuli of definite modulus nesting around the remaining ends of $S$.  Therefore $K = \cx\setminus S$ is a Cantor set of absolute area zero,
and hence $f : S \arrow S$ extends to a polynomial endomorphism of $\chat$
(see e.g. \cite[\S 2.8]{McMullen:book:CDR} and \cite[\S 8D]{Sario:Nakai:book}.)

\bold{The approximation step.}
An arbitrary metrized polynomial-like tree 
$(T,\dist,F)$ in $\cT_\degf$ can be approximated in the geometric topology 
by a sequence $(T_n, \dist_n,F_n)$ of 
escaping trees (Corollary \ref{cor:denseescape}).  
Realize each escaping tree by a polynomial $f_n$.  
The maximal escape rates $M(f_n) =
\dist_n(v_0, J(F_k))$ converge to $\dist(v_0, J(F))$; by 
Proposition \ref{prop:BHproper} these 
polynomials lie in a compact subset of $\MPoly_\degf^*$. 
Pass to a convergent subsequence $[f_n] \arrow [f]$.  By Theorem \ref{thm:cpe}
the tree map $\tau:\MPoly_\degf^*\to\cT_\degf$ is continuous, so $(T,\dist,F)$ is the 
metrized polynomial-like tree associated to $f$.  
This completes the proof of Theorem \ref{thm:realization}.
\qed

\bold{Proof of Theorem \ref{thm:quotmap}.}
Continuity, equivariance, and properness follow from Theorem \ref{thm:cpe}.
Surjectivity is Theorem \ref{thm:realization}.
\qed

\bold{Notes and references.}
For more on the {\em Hurwitz problem} of constructing coverings of
surfaces with specified branching behavior, see e.g. \cite{Edmonds:Kulkarni:Stong},
\cite{Vakil:Hurwitz} and the references therein.
Proposition \ref{prop:hurwitz} above is also covered by
\cite[Thm. 5.2]{Edmonds:Kulkarni:Stong}.


\section{Compactification}
\label{sec:compact}

In this section, we show the projective space of trees $\PT_\degf$
forms a natural boundary for the moduli space of polynomials, 
and that the translation lengths in trees record the limiting
multipliers at periodic points
(Theorems \ref{thm:boundary} and \ref{thm:onemultiplier} of the
Introduction).
As a corollary, we show that the log-multiplier spectra converge to the length 
spectrum of the limiting tree (Theorem \ref{thm:multboundary}). 

\bold{Compactifying moduli space.}  	
Recall that $\tau :\MPoly_\degf^*\to\cT_\degf$ 
assigns to each polynomial with disconnected
Julia set its associated metrized tree map $(T,\dist,F)$.  Projectivizing, we obtain 
a continuous and surjective map to $\proj\cT_\degf$ where the height metric $\dist$
is only determined up to scale.   The map $\tau$ makes 
$\MPoly_\degf \cup \,\proj\cT_\degf$ into a compact
topological space:  every unbounded sequence $[f_n]$ in $\MPoly_\degf$ has a 
subsequence for
which $\tau(f_n)$ converges in $\proj\cT_\degf$.  The following lemma
implies that all points in $\proj\cT_\degf$ arise as limits of polynomials.

\begin{lemma}  \label{lemma:stretch}
The projectivization of $\tau$ to $\proj\cT_\degf$ satisfies
  $$\tau(\MPoly_\degf^*\setminus K) = \tau(\MPoly_\degf^*) = \proj\cT_\degf $$
for every compact $K\subset \MPoly_\degf$.
\end{lemma}

\bold{Proof.} 
The first equality is immediate from Propositions \ref{prop:stretchpoly} 
and \ref{prop:stretch}, and the 
second is the surjectivity of $\tau$ (Theorem \ref{thm:realization}).
\qed

\bold{Proof of Theorem \ref{thm:boundary}.}  
Theorems \ref{thm:cpe} and \ref{thm:projcompact} 
imply that $\proj\cT_\degf$ 
defines a boundary to $\MPoly_\degf$ via the continuous map $\tau$, 
making
$$ \closure{\MPoly}_\degf = \MPoly_\degf \cup \,\proj\cT_\degf$$
into a compact topological space;  Lemma \ref{lemma:stretch}
shows that $\MPoly_\degf$ is dense in $\closure{\MPoly}_\degf$.

To see that iteration $[f]\mapsto [f^n]$ extends continuously to this boundary,
first note that iteration  $i_n: \cT_\degf \to \cT_{\degf^n}$,
defined by $(T,\dist, F)\mapsto (T, \dist, F^n)$, is continuous in the geometric topology
(Lemma \ref{lemma:iteration}),
and $i_n(\tau(f)) = \tau(f^n)$. 
Suppose $[f_k]$ is a sequence in $\MPoly_\degf^*$ converging to 
the normalized tree $(T,\dist, F)\in \del\MPoly_\degf$.  Then $[f_k]$ is unbounded 
in $\MPoly_\degf$, so $M(f_k)\to \infty$ by Proposition \ref{prop:BHproper},
and the normalized trees $(T_k, \dist_k, F_k)$ associated to $f_k$ converge
to $(T,\dist, F)$.  Consequently, 
 $$M(f_k^n) = M(f_k) \to \infty$$
for each $n$, so $[f_k^n]$ is unbounded in $\MPoly_{\degf^n}$, and 
 $$(T_k, \dist_k, F^n_k) \to (T, \dist, F^n)$$
in $\cT_{\degf^n}$.  Thus, we have 
$[f_k^n] \to (T, \dist, F^n)$
in $\closure{\MPoly}_{\degf^n}$.
\qed

\bold{Limits of multipliers.}
We now turn to the proof of Theorem \ref{thm:onemultiplier}.

To formulate this theorem more precisely, it is useful to
introduce the bundle $\MPoly_{\degf,1} \arrow \MPoly_\degf$ of pairs $[f,p]$,
where $f(z)$ is a polynomial of degree $\degf$ and $p$ is a point in 
$\cx$.  For any affine transformation $A(z) = az+b$, we
regard $[f,p]$ and $[AfA^{-1},A(p)]$ as representing the
same point in $\MPoly_{\degf,1}$.  There is a natural continuous map
\begin{displaymath}
	\MPoly_{\degf,1}^* \arrow \cT_{\degf,1}
\end{displaymath}
sending $[f,p]$ to the tree $(T,\dist,F) = \tau(f)$
with the marked point $\pi(p) \mem \Tbar$.
This projection makes $\proj \cT_{\degf,1}$ into a 
into a natural boundary for $\MPoly_{\degf,1}$,
compatible with the compactification of $\MPoly_\degf$ by $\proj\cT_\degf$.

\begin{theorem}
Let $[f_k,z_k] \mem \MPoly_{\degf,1}$ be a sequence of polynomials 
with distinguished periodic points, satisfying $f_k^n(z_k) = z_k$.
Suppose $[f_k,z_k]$ converges to the normalized pointed tree 
$(T,\dist,F,p) \mem \proj \cT_{\degf,1}$.
Then $F^n(p)=p$, and we have:
\begin{equation}
\label{eq:multform}
	L(p,F^n) = \lim_{k\to\infty} \frac{\log^+|(f_k^n)'(z_k)|}{M(f_k)} .
\end{equation}
\end{theorem}

\bold{Proof.}
Let $(T_k,\dist_k, F_k)$ denote the normalized tree maps obtained by
rescaling $\tau(f_k)$ so that $\dist_k(v_0, J(F_k))=1$.  
By Proposition \ref{prop:lengthconv}, the point $p$ is fixed
by $F^n$ and 
  $$L(\pi(z_k), F_k^n) \arrow L(p, F^n)$$
as $k\to\infty$.  Theorem \ref{thm:mult} implies that 
 $$\log^+|(f_k^n)'(z_k)| = M(f_k)L(\pi(z_k), F_k^n) + O(1) .$$
By Proposition \ref{prop:BHproper}, $M(f_k)\to\infty$ as $k \arrow \infty$,  
and equation (\ref{eq:multform}) follows.
\qed

\bold{Spectra.}
Let $f: \cx\to\cx$ be a polynomial of degree $\degf\geq 2$, and let $\cL_n(f)$
be the unordered collection of log-multipliers $\log^+|(f^n)'(p)|$ for 
periodic points $p$ of period $n$ (counted with multiplicity).  
Note that the cardinality of the set $\cL_n(f)$ is exactly $\degf^n$.  The 
{\em log-multiplier spectrum} of $f$ is the list
 $$\cL(f) = \{ \cL_1(f), \cL_2(f), \cL_3(f), \ldots \}.$$

For a metrized polynomial-like tree $(T,\dist,F)$ of degree $\degf$, and 
for each period $n\geq 1$, let $\cL_n(T,\dist,F)$ be the unordered
collection of translation lengths $L(p,F^n)$ of all periodic points $p\in J(F)$
of period $n$ (with 
respect to the modulus metric $\delta$).  
We count critical periodic ends with multiplicity ($= \deg(p,F^n)$), 
so that $\cL_n(T,\dist,F)$ contains exactly $\degf^n$ 
elements.  The {\em length spectrum} of $(T,\dist,F)$ is the list 
  $$\cL(T,\dist,F) = \{ \cL_1(T,\dist,F), \cL_2(T,\dist,F), \cL_3(T,\dist,F), \ldots \}.$$

\begin{theorem} \label{thm:multboundary}
For any sequence $f_k$ in $\Poly_\degf$ such that 
$$[f_k] \arrow (T,\dist,F) \in \del{\MPoly}_\degf,$$
the log-multiplier spectrum $\cL(f_k)$ converges (up to scale) to the 
length spectrum $\cL(T,\dist,F)$.   
\end{theorem}

\bold{Proof.} 
Let $z_k$ be a sequence of fixed points of $f_k^n$ for some iterate $n$,
and let $(T_k,\dist_k,F_k) = \tau(f_k)$.
Pass to a subsequence so that 
the normalized pointed trees $(T_k,\dist_k, F_k, \pi(z_k))$ converge 
in $\cT_{\degf,1}$ to $(T,\dist, F, p)$.  
By Proposition \ref{prop:lengthconv}, this point $p\in\Tbar$ is fixed by $F^n$;
Theorem \ref{thm:onemultiplier} implies that the log-multipliers of $z_k$ 
converge (up to scale) to the translation length of $p$. 

On the other hand, suppose $p\in J(F)$ is fixed by $F^n$ for some $n$. 
By passing to a strongly convergent subsequence, we can always find 
a sequence $p_k\in J(F_k)$ fixed by $F_k^n$ which converges
to $p$.  For every fixed point $p_k$ of $F_k^n$, there exists at least 
one point $z_k$ fixed by $f_k^n$ such that $\pi(z_k)=p_k$.  
\qed



\section{Trees as limits of Riemann surfaces}
\label{sec:GH} 

In this section we recall the 
Gromov-Hausdorff topology on metric spaces, and observe
that it coincides with the geometric topology on the space of trees.
We prove:

\begin{theorem}
\label{thm:GHlimit}
Let $[f_n] \mem \MPoly_\degf^*$ be a sequence of polynomials converging
to the normalized tree $(T,\dist,F)$ in  $\del\MPoly_\degf = \PT_\degf$.
Then the pointed dynamical metric spaces  
$$( \Omega(f_n), c(f_n), |\omega_n|/M(f_n), f_n)$$
converge in the Gromov-Hausdorff topology to 
$(T, v_0, \dist, F)$.
\end{theorem}
An analogous theorem for Kleinian groups appears in \cite{Paulin:rtrees}.

\bold{Pointed metric spaces.}
Suppose $(X,d_X)$ and $(Y, d_Y)$ are metric spaces with basepoints 
$x_0\in X$ and $y_0\in Y$.  
An {\em $\epsilon$-identification} between $(X,x_0,d_X)$ and $(Y,y_0, d_Y)$ 
is a relation ${\cal R} \subset X\times Y$ such that 
\begin{enumerate}
\item	$(x_0,y_0)\in \cal R$
\item 	if $(x,y)$ and $(x',y')$ are in $\cal R$, then 
  $$|d_X(x,x') - d_Y(y,y')| \leq \epsilon,$$
\item	for every $p \in X$, there exists a pair $(x,y)\in \cal R$ such that 
$d_X(p,x)\leq \epsilon$, and 
\item	for every $q\in Y$, there exists a pair $(x,y)\in \cal R$
such that $d_Y(q,y)\leq \epsilon$.  
\end{enumerate}

\bold{Almost-conjugacy.}
A quadruple $(X, x_0, d_X, F)$, where $(X,d_X)$ is a metric space, $x_0\in X$, and
$F: X\to X$, will be called a {\em pointed dynamical metric space}.
An {\em $\epsilon$-identification} between
$(X, x_0, d_X, F)$ and $(Y, y_0, d_Y, G)$ is a relation
${\cal R}\subset X\times Y$ which satisfies conditions (1--4) above and also
\begin{itemize}
\item[5.]		for each $(x,y)\in\cal R$, there exists $(x',y')\in \cal R$ such 
			that $d_X(x', F(x))\leq \epsilon$ and $d_Y(y', G(y))\leq \epsilon$.  
\end{itemize}

\bold{Gromov-Hausdorff topology.}  
Fix a pointed dynamical metric space $(X, x, d_X, F)$.
For $R>0$, let 
 $$X^R  = \{x\in X: d_X(x_0,x) \leq R\},$$
and define
  $$R' = d_X(x_0, F(x_0)).$$
A basis of open sets for the Gromov-Hausdorff topology can be described as follows.
For each $R>R'$ and $\epsilon>0$, the open set $U(R, \epsilon)$ consists of all 
pointed dynamical metric spaces $(Y, y_0, d_Y, G)$ for which there is an 
$\epsilon$-identification between $(X^R\cap F^{-1}(X^R), x_0, d_X, F)$ and
$(Y^R\cap G^{-1}(Y^R), y_0, d_Y, G)$.  

Each element of the space of trees $\cT_\degf$ defines a pointed dynamical
metric space $(T, v_0, \dist, F)$ where $v_0\in T$ is the base and 
$\dist$ is the height metric.
The Gromov-Hausdorff topology on $\cT_\degf$ identifies 
$(T_1, v_0, \dist_1, F_1)$ and $(T_2, v_0, \dist_2, F_2)$ if and 
only if there exists an isometry $i: (T_1,v_0)\arrow (T_2,v_0)$ such that $F_2\compos i
= i\compos F_1$.   

\begin{lemma}
The geometric topology on $\cT_\degf$ coincides with the Gromov-Hausdorff 
topology.
\end{lemma}

\bold{Proof.} 
Suppose $(T_n,\dist_n,F_n)$ converges strongly to $(T,\dist,F)$.  Then for 
each $k>0$, $\epsilon>0$ and all $n>n(k,\epsilon)$, the simplicial isomorphisms between 
$T_n(k)$ and $T_{n+1}(k)$ can be chosen as $\epsilon$-identifications.  
Thus the sequence $(T_n, v_0, \dist_n, F_n)$ converges 
to $(T, v_0, \dist , F)$ in the Gromov-Hausdorff topology. 

Now suppose that $(T_n,\dist_n,F_n)$ does {\em not} converge to $(T,\dist,F)$
in the geometric topology, but the distances $\dist_n(v_0, F_n(v_0))$ remain
bounded away from 0 and $\infty$.  
Then there exists a strongly convergent subsequence which converges 
to $(S, \dist_S, G)\not= (T, \dist, F)$, and therefore the subsequence converges to
$(S, v_0, \dist_S, G)$ in the Gromov-Hausdorff topology.  
The pointed dynamical metric spaces
$(T, v_0, \dist , F)$ and $(S, v_0, \dist_S, G)$ are distinct, and therefore 
the sequence $(T_n, v_0, \dist_n, F_n)$ does not
converge to $(T, v_0, \dist, F)$ in the Gromov-Hausdorff topology.  
\qed

\bold{The modulus metric.}  The Gromov-Hausdorff topology on the space of
metrized polynomial-like trees $(T,v_0, \delta, F)$ with the modulus metric
in place of the height metric
does {\em not} coincide with the geometric topology.  
For example, consider a degree 3 polynomial $f$ with one
critical point escaping and one in a fixed component of $K(f)$.  Suppose that under 
a slight perturbation (in $\MPoly_3$) the second critical point also escapes.  Then
for the perturbed map, all ends of the tree are at infinite distance from the base, 
but for $T$ the distance to the fixed end is finite.  On the other hand,
we know the trees with the height metric vary continuously with $f$ 
in the geometric topology by Theorem \ref{thm:cpe}.

\bold{Pointed trees.}
The geometric topology on the space of pointed trees
$\cT_{\degf,1}$ defined in \S\ref{sec:topology} is equivalent
to the Gromov-Hausdorff topology on $\cT_{\degf,1}$ where the metric 
space is the completion $\Tbar$ and basepoint is  
the marked point $p\in \Tbar$.

\bold{Convergence of the basins of infinity.}
Let $f$ be a polynomial of degree $\degf\geq 2$ with maximal escape
rate $M(f)>0$.   Let $c(f)$ denote any choice of critical point with 
$G(c(f)) = M(f)$.   As in the Introduction,
there is a flat conformal metric (with singularities) $|\omega|$ on 
$\Omega(f)$ defined by $\omega = 2 \partial G$.  
With normalized metric $|\omega|/M(f)$, 
the total length of any level set $G^{-1}(t)$ with $t>0$ 
is equal to $2\pi/M(f)$.  Therefore,
the quotient $\pi:\Omega(f)\arrow T$ defines a $2\pi/M(f)$-identification 
between $(\Omega(f), c(f), |\omega|/M(f), f)$ and $(T, v_0, \dist, F)$.

\bold{Proof of Theorem \ref{thm:GHlimit}.}    
Let $(T_n,\dist_n,F_n)$ be the normalized polynomial-like tree 
associated to $f_n$.
The maximal escape rate $M(f)$ is a proper function 
of $f\in\MPoly_\degf$ by Proposition \ref{prop:BHproper}.  Therefore, 
any Gromov-Hausdorff limit of the pointed
dynamical metric spaces $(\Omega(f_n),c(f_n), |\omega_n|/M(f_n), f_n)$ 
will coincide with the Gromov-Hausdorff limit of the 
trees $(T_n, v_0, \dist_n, F_n)$.  
But the sequence $(T_n, \dist_n,F_n)$ converges to $(T,\dist,F)$
in the geometric topology, so in fact
 $$(\Omega(f_n), c(f_n), |\omega_n|/M(f_n), f_n) \arrow (T, v_0, \dist, F)$$
in the Gromov-Hausdorff topology as $n\arrow \infty$. 
\qed


\section{Families of polynomials over the punctured disk}
\label{sec:families}

Let $\cO$ be the local ring consisting of germs of analytic 
functions $a(t) = \sum_0^\infty a_n t^n$
defined near $t=0$ in $\cx$, and let $K$ be its field of fractions.  The elements of $K$
are given by Laurent series $a(t) = \sum_{-m}^\infty a_n t^n$,
with finite tails and positive radius of convergence.

Let $f \mem K[z]$ be a monic centered polynomial of degree $\degf \ge 2$.
Then $f$ determines a meromorphic family of polynomials
\begin{displaymath}
	f_t(z) = z^\degf + a_2(t) z^{\degf-2} + \cdots + a_\degf(t)
\end{displaymath}
defined near $t=0$.  Conversely, any such family gives an element of $K[z]$.
(We also note that any $f \mem K[z]$ of degree $\degf \ge 2$ is conjugate to a monic
centered polynomial, so there is no loss of generality with this normalization.)

In this section we establish Theorem \ref{thm:algebraic} by showing:

\begin{theorem}
\label{thm:t0}
Either $f_t(z)$ extends holomorphically to $t=0$, or there is a unique
normalized tree such that
\begin{displaymath}
	\lim_{t \arrow 0} \;[f_t] = (T,\dist,F) \mem \PT_\degf .
\end{displaymath}
In the latter case the edges of $T$ have rational length, and hence
the translation lengths of all of its periodic points are also rational.
\end{theorem}

\bold{Coefficient bounds.}  
For a polynomial $f$ of degree $\degf$, let 
$$G^n(z) = \frac{1}{\degf^n} \log^+ |f^n(z)| , $$ 
so that its escape-rate function  is given by $G(z) = \lim G^n(z)$.
Recall that $M(f) = \max\{G(c): f'(c)=0\}$ denotes the maximal escape rate 
of $f$.  
\begin{lemma}  \label{Max}
The coefficients and maximal escape rate of any 
monic and centered polynomial 
$f(z) =  z^\degf + a_2 z^{\degf-2} + \cdots + a_{\degf-1} z + a_0$
are related by:
  $$\max_i |a_i| = O(e^{\degf M(f)})$$
\end{lemma} 

\bold{Proof.} 
Let $h$ be the conformal map that conjugates $z\mapsto z^\degf$ to $f(z)$ near infinity with derivative 1 at $\infty$.  Then $h$ is univalent on the disk $|z | > R = e^{M(f)}$, and the Green's function is given by $G(z) = \log|h^{-1}(z)|$. Moreover,  $h(z ) = z + b_1/z + b_2/z^2 + \cdots$ since $f$ is monic and centered. By the Koebe estimates, 
all critical points of $f$ lie inside the disk $|z | < 4R$.  In particular,
the coefficients of $f'(z)$ are bounded by $O(R^{\degf-1})$.
This gives  the desired bound
(cf.  \cite{Branner:Hubbard:cubicsI}).
\qed

\begin{lemma}  \label{Gn}
For any monic and centered polynomial 
$f(z) = z^\degf + a_2 z^{\degf-2} + \cdots + a_{\degf-1} z + a_\degf$ of degree $\degf\geq 2$, 
we have
  $$|G(z) - G^n(z)| \leq \frac{1}{\degf^n} \log \left(1+\sum|a_i| \right)$$
for all $z\in\cx$.
\end{lemma}

\bold{Proof.} 
By definition,
  $$G^{n+1}(z) - G^n(z) = \frac{1}{\degf^{n+1}} ( \log^+ |f(f^n(z))| - \degf \log^+ |f^n(z)|).$$
For $z$ such that $|f^n(z)|\leq 1$, we obtain
  $$|G^{n+1}(z) - G^n(z)| \leq \frac{1}{\degf^{n+1}} \log 
  \left( 1 + \sum|a_i|\right).$$
For $z$ such that $|f^n(z)| > 1$, we have
 $$|G^{n+1}(z) - G^n(z)| \leq \frac{1}{\degf^{n+1}} \left| \log \frac{|f(f^n(z)|}{|f^n(z)|^\degf} \right|
 	\leq \frac{1}{\degf^{n+1}} \log \left(1 + \sum|a_i|\right).$$
In either case, we obtain by summing over $n$ that 
\begin{displaymath}
  |G(z) - G^n(z)| \leq \frac{1}{\degf^n} \log \left( 1 + \sum_i|a_i| \right).
\end{displaymath}
\qed

\bold{Meromorphic families.}
Now suppose that $f_t(z)$ is a meromorphic family 
of monic and centered polynomials, whose coefficients are holomorphic
on the punctured disk 
\begin{displaymath}
	\Delta^* = \{t \st 0 < |t| < 1 \} .
\end{displaymath}
Suppose that the critical points of $f_t(z)$ can 
be holomorphically parameterized by functions $c_j(t)$, $j=1,2,\ldots,\degf-1$.
We also suppose that the conjugacy classes $[f_t]$ are unbounded in $\MPoly_\degf$
as $t\to 0$.
Thus there exists an $N\geq 1$ such that $\max_i |a_i(t)| \geq C|t|^N > 0$.
We write $G_t$ for the escape-rate function of $f_t$.  

\begin{prop}  \label{prop:rational}
There exist rational numbers $e_j \mem \zed[1/\degf]$
such that 
$$ G_t(c_j(t)) = e_j \log \frac{1}{|t|} + o \left(\log \frac{1}{|t|}\right) $$
as $|t|\to 0$,
and $e_j>0$ for some $j$.
\end{prop}

\bold{Proof.} 
Fix $j$ and consider the meromorphic function
\begin{displaymath}
	f^n_t(c_j(t)) \asyto C_n t^{-p_n}
\end{displaymath}
as $t \arrow 0$.  If $p_n \arrow \infty$ then for all $n$ sufficiently large, 
we have $p_{n+1} = \degf p_n$ and $C_{n+1} = C_n^\degf$ (this holds, for example, 
provided $p_n > \max_i \{\ord_0 a_i(t)\}$).  Then the limit  
 $$e_j := \lim_{n\to\infty}  p_n/\degf^n $$
exists and is non-zero in the ring $\zed[1/\degf]$.  If the degrees 
$p_n$ remain bounded, we set $e_j=0$.  

For any large $n$ we have
$$	G^n_t(c_j(t)) = \frac{1}{\degf^n} \log^+ |f_t^n(c_j(t))|  = \degf^{-n} p_n \log\frac{1}{|t|} + o\left(\log\frac{1}{|t|}\right) $$
as $t \arrow 0$.  On the other hand
$\max_i |a_i(t)| = O(|t|^{-N})$ for some $N$, so by Lemma \ref{Gn}
$$	G_t(c_j(t)) = G^n_t(c_j(t)) + O\left(N\degf^{-n} \log \frac{1}{|t|}\right) $$
which gives the desired asymptotics.

To prove at least one $e_j>0$, use the fact that
$\max |a_i(t)|$ grows like $|t|^{-N} \gg 0$ to conclude by Lemma \ref{Max} that
$$	\max_j G_t(c_j(t)) \ge C \log \frac{1}{|t|} > 0.$$
There are only finitely many critical points, so for some $j$ 
we must have $e_j>0$.  
\qed

\bold{Proof of Theorem \ref{thm:t0}.}
Let $f_t(z)$ be be a meromorphic family of monic
centered polynomials defined near $t=0$.
If the conjugacy classes $[f_t]$ remain bounded in $\MPoly_\degf$ as $t\to 0$, then $f_t$ extends to
a holomorphic family at $t=0$ (equivalently, $f \mem \cO[z]$).

Now assume that $[f_t]$ is unbounded in $\MPoly_\degf$.  
Write $f_t'(z) = df_t/dz$ 
as a product of irreducible polynomials $g_i(z)$ in the ring $K[z]$, where $K$ is the
field of Laurent series.
Then the discriminant of $g_i$ with respect to $z$ is a
meromorphic function $d_i = d_i(t) \mem K$ which is not identically zero.
Since the zeros of $d_i(t)$ are isolated,
after a base change of the form $t \mapsto st$, $0<s \ll 1$, we can assume 
$d_i(t) \neq 0$ for all $t \mem \Delta^*$.  This implies the number of roots 
of $f_t'(z)=0$ is a constant, independent of $t$.
These roots are cyclically permutated
under monodromy around $t=0$, so after a further base change of the form 
$t \mapsto t^n$
they can be labelled by holomorphic functions $c_1(t),\ldots,c_{\degf-1}(t)$.
We are then in the setting of Proposition \ref{prop:rational}.

Let $(T_t, \dist_t, F_t)$ be the normalized tree associated to $f_t$.
Then the heights of all
critical points converge as $t \arrow 0$:  the limiting heights are given 
by the values $e_j/(\max_i e_i)$ of Proposition \ref{prop:rational}.  
For each $k>0$, the number of subtrees $(T(k), \dist, F|T(k))$ with critical 
heights $\{e_j/(\max_i e_i): j = 1, \ldots, \degf-1\}$ is finite, and 
therefore the family
$(T_t, \dist_t, F_t)$ must converge by 
the continuity of the tree map.  

Proposition \ref{prop:rational} also implies that the heights $e_j/(\max_i e_i)$
are rational, so all edge lengths in the normalized limiting tree are rational.  
\qed

\bold{Remark.}
For more on the dynamics of rational maps over local rings and fields,
see e.g. 
\cite{Benedetto:local},
\cite{Kiwi:trees},
\cite{Rivera-Letelier:thesis}, \cite{Rivera-Letelier:hype}
and \cite[Ch. 8]{Baker:Rumely:book}.


\section{Cubic polynomials}
\label{sec:cubics}

In this section we discuss the topology and combinatorics of $\PT_\degf$ in the
case $\degf=3$.  In particular we prove:

\begin{theorem}
\label{thm:treeoftrees}
The space $\PT_3$ is the completion of
an infinite simplicial tree $\PT_3(2) \cup \PT_3(1,1)$ by its space of ends $\PT_3(1)$.
\end{theorem}
We also describe an algorithm to construct a concrete model for $\PT_3$.

\bold{Strata.}
As mentioned in the Introduction, the space of cubic trees decomposes into strata
\begin{displaymath}
	\PT_3 = \PT_3(2) \disjunion \PT_3(1,1) \disjunion \PT_3(1).
\end{displaymath}
In the top-dimensional stratum $\PT_3(1,1)$ there are two escaping critical points
with disjoint forward orbits; in $\PT_3(2)$ there are two escaping critical points
whose forward orbits collide; and in $\PT_3(1)$ there is only one escaping critical point. 

In each stratum, the simplicial structure of the tree is locally constant;
only the metric varies.
Thus the components of $\PT(1)$ and $\PT(2)$ are points, and those of
$\PT(1,1)$ are open intervals.

\bold{Levels.}
We define the {\em level} of a cubic tree by
\begin{displaymath}
	\ell(T,\dist,F) = \frac{\log H(v_0)/H(c_0)}{\log 3} ,
\end{displaymath}
where $v_0 \mem T$ is the fastest escaping critical point (as usual),
$c_0 \mem \Tbar$ is the remaining critical point, and $H : \Tbar \arrow [0,\infty)$ is the
height function.

The level gives a continuous map $\ell : \PT_3 \arrow [0,\infty]$ such that
\begin{enumerate}
	\item
	$\ell = \infty$ iff $c_0 \mem J(F)$ iff $(T,\dist,F) \mem \PT_3(1)$;
	\item
	$\ell$ is an integer iff $(T,\dist,F) \mem \PT_3(2)$, in which case $F^\ell(c_0)=v_0$; and
	\item
	$\ell$ maps each component of $\PT_3(1,1)$ homeomorphically to an
	open interval of the form $(n,n+1)$ where $n \ge 0$ is an integer.
\end{enumerate}

Recall that $T(k) \subset T$ denotes the finite subtree spanned by the vertices
with combinatorial height $-kN(F) \le h(v) \le kN(F)$.
If $T$ is normalized so $H(v_0)=1$, then $T(k) = H^{-1}([3^{-k},3^k])$.
By Theorem \ref{thm:subtree}, any tree in $\PT_3$ at level $\ell \le k$ is determined
by its level and by its combinatorial dynamics on the finite subtree $T(k)$.
Since there are only finitely many
possibilities for the combinatorial dynamics, we have:

\begin{prop}
The strata $\PT_3(1,1)$ and $\PT_3(2)$ have only finitely many components below
any given level $\ell$.
\end{prop}

\bold{Proof of Theorem \ref{thm:treeoftrees}.}
We begin by showing that $P = \PT_3(2) \cup \PT_3(1,1)$ is a simplicial tree.

The preceding Proposition shows that $P$ is
a simplicial 1-complex, and that $\ell : P \arrow [0,\infty)$ is a proper
map taking integral values exactly on the $0$-cells of $P$.
As remarked earlier, there is a unique point where $\ell=0$.

By Theorem \ref{thm:subtree}, any tree at level $\ell > 0$ is a limit
of trees with levels $\ell_n < \ell$.  Moreover, if
$(T_n,\dist_n,F_n) \arrow (T,\dist,F) \mem \PT_3(2)$ are normalized trees converging from below,
then by the definition of convergence, there is eventually an isometric embedding of the
part $U_n$ of $T_n$ at height $\ge 3^{-\ell_n}$ into $T$, respecting the dynamics.
By Theorem \ref{thm:subtree} again, this means the combinatorial dynamics of
$F_n|T_n$ is eventually determined by $(T,F)$.  Consequently all $(T_n,\dist_n,F_n)$
lie in a single component of $\PT_3(1,1)$ for $n$ sufficiently large.

This shows every $0$-cell of $P$ is the limit from below of a unique $1$-cell,
and thus $P$ is an infinite simplicial tree.  

By density of escaping trees (Corollary \ref{cor:denseescape}), $\closure{P} = \PT_3$.
Given a normalized tree $(T,\dist,F)$ in $\PT_3(1)$, the
dynamics on the subtree of height $\ge 3^{-\ell}$ determines the unique
vertex in $\PT_3(2)$ lying above it at level $\ell$, 
and hence the stratum
$\PT_3(1)$ corresponds canonically to the set of ends of $P$.
\qed

\bold{Encoding trees.}
Next we define a finite sequence of integers encoding each vertex of $\PT_3$.

Let $(T,\dist,F)$ be a tree at level $\ell$ in $\PT_3(2)$, with critical 
points $v_0$ and $c_0$.  As usual we let $v_1, v_2, \ldots$ denote the unique vertices
with combinatorial height $h(v_i) = i$; since $N(F)=1$, we have $v_i = F^i(v_0)$.
The {\em spine} of $T$ is the unique sequence of consecutive vertices 
$(s_0,s_1,\ldots,s_\ell)$ with $s_0=v_0$ and $s_\ell = c_0$.

Let $S \subset T$ be the union of the forward orbits of the vertices $s_i$ in the spine.
Define $g : S \arrow \{0,1,2,\ldots\}$ by setting $g(v)$ equal to the {\em least}
$i$ such that $v$ is in the forward orbit of $s_i$.  For example, $g(s_i) = i$
and $g(v_i)=0$.

Let $N \subset S$ be the set of vertices at combinatorial distance $\le 1$
from the unique path $(s_\ell,s_{\ell-1},\ldots,s_1,s_0=v_0,v_1,v_2,\ldots)$ joining $c_0$
to infinity.  Let $F_N : N \arrow N$ be the first return map; in other words,
$F_N(v) = F^j(v)$ for the least $j>0$ such that $F^j(v) \mem N$.
Then for $1 \le i \le \ell$ we define:
\begin{itemize}
	\item
	the {\em lifetime} $k(i)$ of $s_i$ to be the number of vertices in $N$ with $g(v)=i$; and
	\item
	the {\em terminus} $t(i)$ of $s_i$ to be the value of $g(F_N^j(s_i))$ for the least $j>0$
	such that $g(F_N^j(s_i)) \neq i$.
\end{itemize}
The {\em code} for the tree $(T,\dist,F)$ is the finite sequence $(k(i),t(i))$, $i=1,\ldots,\ell$.

For a tree with only one escaping critical point, the code is the infinite sequence defined 
using the unique path $(s_0,s_1,s_2,\ldots)$ joining $v_0$ to $c_0 \mem J(F)$.

\begin{theorem}
Any tree $(T,\dist,F) \mem \PT_3(2) \cup \PT_3(1)$ is uniquely determined by its code.
\end{theorem}

\bold{Proof.}
The data $(k(i),t(i))$ allow one to inductively reconstruct the first-return map $F_N : N\arrow N$
and the combinatorial height function $h : N \arrow N$.  
This data in turn determines $F$ on the forward-invariant subtree spanned by $S$,
which determines $(T,\dist,F)$ by Theorem \ref{thm:subtree}.
\qed

\makefig{Codes for all 42 cubic trees at level 6.}{fig:codes}{
{\scriptsize
\begin{tabular}{cccccccccccc}
  1 & 0 & 1 & 1 & 1 & 2 & 1 & 3 & 1 & 4 & 1 & 5\\
  1 & 0 & 1 & 1 & 1 & 2 & 1 & 3 & 1 & 4 & 6 & 0\\
     \stroke{8}{4}
  1 & 0 & 1 & 1 & 1 & 2 & 1 & 3 & 5 & 0 & 1 & 1\\
  1 & 0 & 1 & 1 & 1 & 2 & 1 & 3 & 5 & 0 & 1 & 5\\
     \stroke{6}{6}
  1 & 0 & 1 & 1 & 1 & 2 & 4 & 0 & 1 & 1 & 1 & 2\\
  1 & 0 & 1 & 1 & 1 & 2 & 4 & 0 & 1 & 1 & 1 & 4\\
  1 & 0 & 1 & 1 & 1 & 2 & 4 & 0 & 1 & 1 & 2 & 4\\
     \stroke{8}{4}
  1 & 0 & 1 & 1 & 1 & 2 & 4 & 0 & 1 & 4 & 1 & 1\\
  1 & 0 & 1 & 1 & 1 & 2 & 4 & 0 & 1 & 4 & 1 & 4\\
     \stroke{4}{8}
  1 & 0 & 1 & 1 & 3 & 0 & 1 & 1 & 1 & 2 & 1 & 3\\
  1 & 0 & 1 & 1 & 3 & 0 & 1 & 1 & 1 & 2 & 2 & 2\\
  1 & 0 & 1 & 1 & 3 & 0 & 1 & 1 & 1 & 2 & 3 & 3\\
     \stroke{8}{4}
  1 & 0 & 1 & 1 & 3 & 0 & 1 & 1 & 1 & 3 & 1 & 1\\
  1 & 0 & 1 & 1 & 3 & 0 & 1 & 1 & 1 & 3 & 1 & 3\\
     \stroke{8}{4}
  1 & 0 & 1 & 1 & 3 & 0 & 1 & 1 & 2 & 3 & 1 & 1\\
  1 & 0 & 1 & 1 & 3 & 0 & 1 & 1 & 2 & 3 & 1 & 3\\
     \stroke{6}{6}
  1 & 0 & 1 & 1 & 3 & 0 & 1 & 3 & 1 & 1 & 1 & 2\\
  1 & 0 & 1 & 1 & 3 & 0 & 1 & 3 & 1 & 1 & 1 & 3\\
  1 & 0 & 1 & 1 & 3 & 0 & 1 & 3 & 1 & 1 & 2 & 3\\
     \stroke{8}{4}
  1 & 0 & 1 & 1 & 3 & 0 & 1 & 3 & 1 & 3 & 1 & 1\\
  1 & 0 & 1 & 1 & 3 & 0 & 1 & 3 & 1 & 3 & 1 & 3\\
     \stroke{2}{10}
  1 & 0 & 2 & 0 & 1 & 1 & 1 & 2 & 1 & 3 & 1 & 4\\
  1 & 0 & 2 & 0 & 1 & 1 & 1 & 2 & 1 & 3 & 3 & 1\\
     \stroke{8}{4}
  1 & 0 & 2 & 0 & 1 & 1 & 1 & 2 & 2 & 2 & 1 & 1\\
  1 & 0 & 2 & 0 & 1 & 1 & 1 & 2 & 2 & 2 & 1 & 2\\
     \stroke{6}{6}
  1 & 0 & 2 & 0 & 1 & 1 & 2 & 1 & 1 & 2 & 1 & 3\\
  1 & 0 & 2 & 0 & 1 & 1 & 2 & 1 & 1 & 2 & 2 & 2\\
     \stroke{8}{4}
  1 & 0 & 2 & 0 & 1 & 1 & 2 & 1 & 1 & 4 & 1 & 2\\
  1 & 0 & 2 & 0 & 1 & 1 & 2 & 1 & 1 & 4 & 1 & 4\\
  1 & 0 & 2 & 0 & 1 & 1 & 2 & 1 & 1 & 4 & 2 & 1\\
     \stroke{8}{4}
  1 & 0 & 2 & 0 & 1 & 1 & 2 & 1 & 2 & 1 & 1 & 2\\
  1 & 0 & 2 & 0 & 1 & 1 & 2 & 1 & 2 & 1 & 1 & 4\\
  1 & 0 & 2 & 0 & 1 & 1 & 2 & 1 & 2 & 1 & 1 & 5\\
     \stroke{4}{8}
  1 & 0 & 2 & 0 & 1 & 2 & 1 & 1 & 1 & 2 & 1 & 3\\
  1 & 0 & 2 & 0 & 1 & 2 & 1 & 1 & 1 & 2 & 2 & 1\\
     \stroke{8}{4}
  1 & 0 & 2 & 0 & 1 & 2 & 1 & 1 & 2 & 1 & 1 & 2\\
  1 & 0 & 2 & 0 & 1 & 2 & 1 & 1 & 2 & 1 & 1 & 5\\
  1 & 0 & 2 & 0 & 1 & 2 & 1 & 1 & 2 & 1 & 2 & 1\\
     \stroke{6}{6}
  1 & 0 & 2 & 0 & 1 & 2 & 1 & 2 & 1 & 1 & 1 & 2\\
  1 & 0 & 2 & 0 & 1 & 2 & 1 & 2 & 1 & 1 & 2 & 1\\
     \stroke{8}{4}
  1 & 0 & 2 & 0 & 1 & 2 & 1 & 2 & 1 & 2 & 1 & 1\\
  1 & 0 & 2 & 0 & 1 & 2 & 1 & 2 & 1 & 2 & 1 & 2
\end{tabular}
}
\medskip
}

\bold{Examples.}
For simplicity we will present codes in the folded form 
\begin{displaymath}
	(k(1),t(1),k(2),t(2),\ldots).
\end{displaymath}
The empty sequence and the sequence $(1,0)$ encode the unique trees at 
levels $\ell=0$ and $1$.
The code for the quotient tree of $f(z)=cz^2+z^3$, $c \gg 0$,
is the infinite sequence $(1,0,1,1,1,2,1,3,1,4,\ldots)$.

The codes for the 42 cubic trees at level $\ell=6$ in $\PT_3$ are shown
in Table \ref{fig:codes}.  
A list of codes for the vertices at level $6-i$ is obtained by
dropping the last $2i$ columns of the table.
Since the table is sorted, it also depicts the shape of the tree $\PT_3$ itself.
We have added horizontal lines to make the tree structure more visible.

\bold{Tableaux and $\tau$-functions.}
The Yoccoz $\tau$-function 
is defined for $0<i \le \ell$ by 
\begin{displaymath}
	\tau(i) = \max \{j \st s_j = F^k(s_i) \;\text{for some $k>0$}\}
\end{displaymath}
(see e.g. \cite{Hubbard:local:connectivity}, \cite[Problem 1e]{Milnor:local:connectivity}).
The $\tau$-function can be computed inductively
from the code for a tree; namely $\tau(1)=0$ and 
\begin{displaymath}
	\tau(i) = \begin{cases}
		t(i) & \text{if $t(i) = \tau^{k(i)}(i-1)+1$,}\\
		\tau(t(i)) & \text{otherwise,}
		\end{cases}
\end{displaymath}
where we adopt the convention that $\tau(0)=-1$.

On the other hand, the $\tau$-function does not determine the tree.
Already at level $\ell=5$, the trees with codes 
\begin{displaymath}
	(1,0,1,1,3,0,1,1,1,3) \;\;\;\text{and}\;\;\; (1,0,1,1,3,0,1,1,2,3) 
\end{displaymath}
have the same $\tau$-function, namely $\tau(1,2,3,4,5)=(0,1,0,1,0)$.
It is known that the $\tau$-function carries the same information as
the critical tableau introduced in \cite{Branner:Hubbard:cubicsII},
so trees also carry more information than tableaux.

\maketab{The number of vertices in $\PT_3$ at level $\ell$.}{tab:Nl}{
\begin{tabular}{|c|ccccccccc|}\hline
$\ell$	 & 0 & 1 & 2 & 3 & 4 &  5 &  6 &   7 &   8 \\ \hline
$N_\ell$ & 1 & 1 & 2 & 4 & 8 & 18 & 42 & 103 & 260 \\ \hline \hline
$\ell$	 &   9 &   10 &   11 &    12 &    13 &    14 &     15 &     16 &      17 \\ \hline
$N_\ell$ & 670 & 1753 & 4644 & 12433 & 33581 & 91399 & 250452 & 690429 & 1913501 
\\ \hline
\end{tabular}
}

\bold{Growth of $\PT_3$.}
It is straightforward to automate the inductive enumeration of tree codes,
by keeping track of $F_N|N$.
The resulting vertex counts for $\ell \le 17$ are given in
Table \ref{tab:Nl}.

\bold{Question.}
Does $\lim_{\ell \arrow \infty} \;(1/\ell) \log N_\ell$ exist?  Is it equal to $\log 3$?

\bold{Notes.}
We remark that an additional tableau rule is needed in the realization theorem 
\cite[Thm. 4.1]{Branner:Hubbard:cubicsII}.
In terms of the $\tau$-function, the extra condition is that
if $(\tau(i),\tau(i+1)) = (0,0)$ for some $i \ge 1$,
then $(\tau(j),\tau(j+1)) \neq (i,0)$ for all $j \ge 0$.
For example, the tableau associated to $\tau(1,2,3,4)=(0,0,1,0)$ satisfies
the conditions (Ta), (Tb) and (Tc) of \cite[\S 4]{Branner:Hubbard:cubicsII}
but cannot be realized by a cubic polynomial.
See also \cite{Harris:tableau} and \cite[Remark 4.6]{Kiwi:trees}.

\bibliographystyle{math}
\bibliography{math}

\begin{thebibliography}{FLM}

\bibitem[BR]{Baker:Rumely:book}
M.~Baker and R.~Rumely.
\newblock {Potential Theory on the Berkovich Projective Line}.
\newblock {\em In preparation}.

\bibitem[Be]{Benedetto:local}
R.~Benedetto.
\newblock {Reduction, dynamics, and Julia sets of rational functions}.
\newblock {\em J. Number Theory} {\bf 86}(2001), 175--195.

\bibitem[BH1]{Branner:Hubbard:cubicsI}
B.~Branner and J.~H. Hubbard.
\newblock {The iteration of cubic polynomials, Part I: The global topology of
  parameter space}.
\newblock {\em Acta Math.} {\bf 160}(1988), 143--206.

\bibitem[BH2]{Branner:Hubbard:cubicsII}
B.~Branner and J.~H. Hubbard.
\newblock {The iteration of cubic polynomials, Part II: Patterns and
  parapatterns}.
\newblock {\em Acta Math.} {\bf 169}(1992), 229--325.

\bibitem[Bro]{Brolin:measure}
H.~Brolin.
\newblock {Invariant sets under iteration of rational functions}.
\newblock {\em Ark. Math.} {\bf 6}(1965), 103--144.

\bibitem[EKS]{Edmonds:Kulkarni:Stong}
A.~Edmonds, R.~Kulkarni, and R.~E. Stong.
\newblock {Realizability of branched coverings of surfaces}.
\newblock {\em Trans. Amer. Math. Soc.} {\bf 282}(1984), 773--790.

\bibitem[Em1]{Emerson:trees}
N.~D. Emerson.
\newblock {Dynamics of polynomials with disconnected Julia sets}.
\newblock {\em Discrete Contin. Dyn. Syst.} {\bf 9}(2003), 801--834.

\bibitem[Em2]{Emerson:harmonic}
N.~D. Emerson.
\newblock {Brownian motion, random walks on trees, and harmonic measure on
  polynomial Julia sets}.
\newblock {\em Preprint, 9/2006}.

\bibitem[FLM]{FLM:entropy}
A.~Freire, A.~Lopes, and R.~Ma{\~n}\'e.
\newblock {An invariant measure for rational maps}.
\newblock {\em Bol. Soc. Brasil. Mat.} {\bf 14}(1983), 45--62.

\bibitem[Gr1]{Gromov:book:var}
M.~Gromov.
\newblock {\em Structures m\'{e}triques pour les vari\'{e}t\'{e}s
  riemanniennes}.
\newblock CEDIC, Textes math\'{e}matiques, 1981.

\bibitem[Gr2]{Gromov:entropy}
M.~Gromov.
\newblock {On the entropy of holomorphic maps}.
\newblock {\em Enseign. Math.} {\bf 49}(2003), 217--235.

\bibitem[Ha]{Harris:tableau}
David~M. Harris.
\newblock {Turning curves for critically recurrent cubic polynomials}.
\newblock {\em Nonlinearity} {\bf 12}(1999), 411--418.

\bibitem[Hub]{Hubbard:local:connectivity}
J.~H. Hubbard.
\newblock {Local connectivity of Julia sets and bifurcation loci: three
  theorems of J.-C. Yoccoz}.
\newblock In L.~R. Goldberg and A.~V. Phillips, editors, {\em Topological
  Methods in Modern Mathematics}, pages 467--511. Publish or Perish, Inc.,
  1993.

\bibitem[Ki]{Kiwi:trees}
J.~Kiwi.
\newblock {Puiseux series of polynomial dynamics and iteration of complex cubic
  polynomials}.
\newblock {\em Ann. Inst. Fourier} {\bf 56}(2006), 1337--1404.

\bibitem[KS]{Kozlovski:vanStrien}
O.~Kozlovski and S.~van Strien.
\newblock {Local connectivity and quasi-conformal rigidity of
  non-renormalizable polynomials}.
\newblock {\em Preprint, 2006}.

\bibitem[LV]{Lehto:Virtanen:book}
O.~Lehto and K.~J. Virtanen.
\newblock {\em Quasiconformal Mappings in the Plane}.
\newblock Springer-Verlag, 1973.

\bibitem[Ly]{Lyubich:entropy}
M.~Lyubich.
\newblock {Entropy properties of rational endomorphisms of the Riemann sphere}.
\newblock {\em Ergodic Theory Dynam. Systems} {\bf 3}(1983), 351--385.

\bibitem[Mc1]{McMullen:aut}
C.~McMullen.
\newblock {Automorphisms of rational maps}.
\newblock In {\em Holomorphic Functions and Moduli I}, pages 31--60.
  Springer-Verlag, 1988.

\bibitem[Mc2]{McMullen:book:CDR}
C.~McMullen.
\newblock {\em Complex Dynamics and Renormalization}, volume 135 of {\em Annals
  of Math. Studies}.
\newblock Princeton University Press, 1994.

\bibitem[Mc3]{McMullen:curdev}
C.~McMullen.
\newblock {The classification of conformal dynamical systems}.
\newblock In {\em Current Developments in Mathematics, 1995}, pages 323--360.
  International Press, 1995.

\bibitem[Mc4]{McMullen:rtrees}
C.~McMullen.
\newblock {Ribbon $\mathbb R$-trees and holomorphic dynamics on the unit disk}.
\newblock {\em Preprint, 11/2007}.

\bibitem[Mil]{Milnor:local:connectivity}
J.~Milnor.
\newblock {Local connectivity of Julia sets: expository lectures}.
\newblock In Tan Lei, editor, {\em The Mandelbrot Set, Theme and Variations},
  pages 67--116. Cambridge University Press, 2000.

\bibitem[Mor]{Morgan:icm86}
J.~Morgan.
\newblock {Trees and hyperbolic geometry}.
\newblock In {\em Proceedings of the International Congress of Mathematicians
  (Berkeley, 1986)}, pages 590--597. Amer. Math. Soc., 1987.

\bibitem[MS1]{Morgan:Shalen:I}
J.~Morgan and P.~Shalen.
\newblock {Valuations, trees, and degenerations of hyperbolic structures, I}.
\newblock {\em Annals of Math.} {\bf 120}(1984), 401--476.

\bibitem[MS2]{Morgan:Shalen:survey}
J.~Morgan and P.~Shalen.
\newblock {An introduction to compactifying spaces of hyperbolic structures by
  actions on trees}.
\newblock In {\em Geometry and Topology}, volume 1167 of {\em Lecture Notes in
  Mathematics}, pages 228--240. Springer-Verlag, 1985.

\bibitem[Ot]{Otal:book:fibered}
J.-P. Otal.
\newblock {\em Le th\'eor\`eme d'hyperbolisation pour les vari\'et\'es
  fibr\'ees de dimension trois}.
\newblock Ast\'erisque, vol. 235, 1996.

\bibitem[Par]{Parry:book:entropy}
W.~Parry.
\newblock {\em Entropy and generators in ergodic theory}.
\newblock W. A. Benjamin, Inc., 1969.

\bibitem[Pau]{Paulin:rtrees}
F.~Paulin.
\newblock {Topologie de Gromov \'equivariante, structures hyperboliques et
  arbres r\'eels}.
\newblock {\em Invent. math.} {\bf 94}(1988), 53--80.

\bibitem[Pil]{Pilgrim:dessins}
K.~Pilgrim.
\newblock {Dessins d'enfants and Hubbard trees}.
\newblock {\em Ann. Sci. \'Ec. Norm. Sup.} {\bf 33}(2000), 671--693.

\bibitem[PS]{Przytycki:Skrzypczak}
F.~Przytycki and J.~Skrzypczak.
\newblock {Convergence and pre-images of limit points for coding trees for
  iterations of holomorphic maps}.
\newblock {\em Math. Ann.} {\bf 290}(1991), 425--440.

\bibitem[QY]{Qiu:Yin}
W.~Qiu and Y.~Yin.
\newblock {Proof of the Branner-Hubbard conjecture on Cantor Julia sets}.
\newblock {\em Preprint, 2006}.

\bibitem[Ri1]{Rivera-Letelier:thesis}
J.~Rivera-Letelier.
\newblock {Dynamique des fonctions rationnelles sur des corps locaux}.
\newblock In {\em Geometric methods in dynamics. II.}, pages 147--230.
  Ast\'erisque, vol. 287, 2003.

\bibitem[Ri2]{Rivera-Letelier:hype}
J.~Rivera-Letelier.
\newblock {Points p\'eriodiques des fonctions rationnelles dans l'espace
  hyperbolique $p$-adique}.
\newblock {\em Comment. Math. Helv.} {\bf 80}(2005), 593--629.

\bibitem[SN]{Sario:Nakai:book}
L.~Sario and M.~Nakai.
\newblock {\em Classification Theory of Riemann Surfaces}.
\newblock Springer-Verlag, 1970.

\bibitem[Shi]{Shishikura:rings}
M.~Shishikura.
\newblock {Trees associated with the configuration of Herman rings}.
\newblock {\em Ergod. Th. \& Dynam. Sys.} {\bf 9}(1989), 543--560.

\bibitem[Sp]{Springer:book:RS}
G.~Springer.
\newblock {\em Riemann Surfaces}.
\newblock Chelsea Publishing Co., 1981.

\bibitem[Va]{Vakil:Hurwitz}
R.~Vakil.
\newblock {Genus 0 and 1 Hurwitz numbers: recursions, formulas, and
  graph-theoretic interpretations}.
\newblock {\em Trans. Amer. Math. Soc.} {\bf 353}(2001), 4025--4038.

\end{thebibliography}

\bigskip
{
\sc 
Mathematics Department,
University of Illinois at Chicago,
851 S Morgan (M/C 249), Chicago, IL  60607-7045.

\bigskip
Mathematics Department, Harvard University, 1 Oxford St,
Cambridge, MA 02138
}

\end{document}